\documentclass[a4paper,10pt,draft]{article}

\usepackage{babel,amsmath,amssymb,amsbsy,amsfonts,latexsym,amsthm}

\newtheorem{Def}{Definition}[section]

\newtheorem{teo}[Def]{Theorem}

\newtheorem{prop}[Def]{Proposition}

\newtheorem{lem}[Def]{Lemma}

\newtheorem{oss}[Def]{Remark}

\newenvironment{Dim}[0]{{\bf Proof}}{}

\newenvironment{Dimo}[0]{{\bf{Proof of}}}

\def\o{\overline}

\def\w{\widehat}

\def\*{\star}

\def\t{\tilde}

\def\ra{\rangle}

\def\la{\langle}

\def\R{{\mathbb R}}

\def\E{{\mathbb E}}

\def\L1{\Lambda_B^1([0,T])}

\def\L2{\Lambda_B^2([0,T])}

\def\l{\lambda}

\def\L{\Lambda}

\def\M1{M^1_B([0,T])}

\def\M2{M^2_B([0,T])}

\def\W{{\cal W}}

\def\M{{\cal M}}

\def\C{{\cal C}}

\def\F{{\cal F}}

\def\s{\sigma}

\def\p{\prime}

\def\o{\overline}

\def\t{\tilde}

\def\P{{\cal P}}

\def\C{{\cal C}}

\def\S{{\cal S}}

\def\la{\langle}

\def\ra{\rangle}

\def\wx{\widehat{x}}

\def\wy{\widehat{y}}

\def\wz{\widehat{z}}

\def\wt{\widehat{t}}

\def\var{\varepsilon}

\def\acc{\'}

\def\R{{\mathbb R}}

\def\W{{\mathbb W}}

\def\E{{\mathbb  E}}

\newcommand{\cvd}{\begin{flushright}

\rule[5pt]{5pt}{5pt}

\end{flushright}}

\date{}

\begin{document}

\title{Regularity results for a class of Semilinear Parabolic Degenerate Equations and Applications}

\maketitle

\begin{center}

\author{Marco Papi\\

\vspace{10pt}

\begin{scriptsize}
Istituto per le Applicazioni del Calcolo "M.Picone", V.le del Policlinico 137, I-00161 Roma (Italy) and\\ 
Dipartimento di Matematica di Roma "Tor Vergata", Via della Ricerca Scientifica, 00133, (Italy) 
\end{scriptsize}} 

\end{center}

\vspace*{10pt}

\begin{center}

{\bf Abstract}

\end{center}

\begin{quote}

\begin{small} 
We prove some regularity results for viscosity solutions to strongly degenerate parabolic semilinear 
problems. These results apply to a specific model used for pricing Mortgage-Backed Securities and allow
a complete justification of use of the classical Ito's formula.
\end{small}

\end{quote}

\vspace*{10pt}
\begin{quote}
\begin{small}
{\bf AMS subject classifications:} 35K55, 35K65 , 35B65, 60H30, 91B70. 
\end{small}
\end{quote}

\section{Introduction}
In this paper we investigate the second order regularity of bounded viscosity solutions to the following semilinear parabolic equation of degenerate type
\begin{eqnarray}\label{EQ}
\partial_t u+H(x,t,u,\nabla u,\nabla^2 u)=0,
\end{eqnarray}
where the hamiltonian function is defined by the expression,
\begin{eqnarray}\label{H}
H(x,t,u,p,X)&=&-\frac{1}{2}tr(\s\s^{\top}(t)X)+\la\mu(x,t), p\ra+\l(u)|\s^{\top}p|^2\nonumber\\
&&+\eta(u)\la\s^{\top}(t)p,w(x,t)\ra+f(x,t,u),
\end{eqnarray}
for every $(x,t)\in\R^N\times(0,T)$, $a<u<b$, $p\in\R^N$, $X\in\S^N$, where $\S^N$ is the space of $N\times N$ symmetric matrices endowed with the usual ordering. In (\ref{H}), $tr$, $|\cdot|$ and $\la,\ra$ denote the 
trace of a square matrix, the Euclidean norm and inner product, respectively. Moreover  if $\M_{N\times d}(\R)$ is the space of  real $N\times d$
matrices, with $N\geq d$, then we assume that, $\mu:\R^N\times[0,T)\rightarrow \R^N$, $\s:[0,T)\rightarrow \M_{N\times d}(\R)$, $w:\R^N\times[0,T)\rightarrow \R^d$, and $f:\R^N\times[0,T)\times(a,b)\rightarrow \R$, are continuous functions.\\
Actually our main motivation comes from the following semilinear equation

\begin{eqnarray}\label{DM1}
\!\partial_t U\!-\!\frac{1}{2}tr(\s\s^{\top}\nabla^2 U)\!-\!\la \mu,\nabla U \ra\!+\!\rho
\frac{|\s^{\top}\nabla U|^2}{U+h+\xi(t)}
+r(U\!+\!h)\!-\!\tau h=0,
\end{eqnarray}
in $\R^N\times[0,T),\;\rho>0,\;\tau,\;T>0$, where $\xi=\xi(t)$ and $h=h(x,t)$ are regular functions of their variables. 
Equation (\ref{DM1}) has been proposed in \cite{37} as a differential model used for pricing some widely traded American financial instruments, the Mortgage-Backed 
Securities ($MBS$); in particular the model was derived following the outline of X. Gabaix in \cite{19}.\\
Let $U$ be a viscosity solution of (\ref{DM1}), then setting $u=U+h+\xi$, $u$ solves a differential equation of type (\ref{EQ}). The arbitrage pricing principle applies to financial instruments whose cash flows are related to the values of some economic factors, such as the interest rates ($r=r(t)$). Using that principle the value a these securities can
be expressed as a conditional expectation over the probability space ($\Omega$) of the underlying factors
that determine the instrument's price ($V_t$) with respect to 
a particular measure ($Q$) defined over that space. As a consequence the knowledge of second order regularity properties of $u$, and therefore of $U$, plays a fundamental role in order to
close in a rigorous way the financial argument used for deriving equation (\ref{DM1}). If we have sufficient conditions about $U$, for applying the classical Ito's rule (see for istance
\cite{9}, or \cite{10}), then, as we will show in section \ref{Conclusions}, there exists a probability measure called the {\em Equivalent Martingale Measure}, $Q$, which has the 
required empiric and statistic representation proposed in \cite{19}, such that the following probabilistic equality holds:
\begin{eqnarray}\label{GABAIX}
U(X_t,T-t)=\E_t^{Q}\Big[\int_t^T \big(\tau-r(T-s)\big) e^{-\int_t^s r(T-\kappa)d\kappa} h(X_s,T-s) ds\Big],\;\;\mbox{a.s}.\end{eqnarray}
\begin{eqnarray}\label{GABAIX.1}
d X_t = \mu(X_t,T-t)dt+\s(T-t)dW_t,\;\;\;0<t\leq T.
\end{eqnarray}
There $\E_t^Q[\cdot]$ denotes the conditional mean value up to the time $t$, taken with respect to $Q$, and $\xi(t)=e^{\int_0^t r(s)ds}$, $U(\cdot,0)\equiv 0$, $h\geq 0$, $h(\cdot,0)\equiv 0$, $r$ is a deterministic function which represents the interest rate variable. 
The process $W_t$ is a $d$-dimensional standard Brownian Motion over the probability space $\Omega$, while $X_t$ is a $N$-dimensional Ito process which describes
the factors which affect the value of a Mortgage-Backed security. In this model $h$ models the remaining principal during the life 
of the pool of mortgages (see \cite{19}, \cite{20}).\\
The equation (\ref{GABAIX}) represents the conclusive statement of our works about the pricing equation introduced in \cite{37} and then studied from the point of view of the existence and uniqueness. Actually that formula gives the formalization of the existence of the market price of risk proposed by X. Gabaix in \cite{19}, which describes the martingale measure such that the value of the security can be expressed as the solution of the Hamilton Jacobi equation (\ref{DM1}).\\
Our main result concerns the semiconvexity/concavity property, at a fixed time $t\in[0,T)$, of the viscosity solution $u$ of the equation (\ref{EQ}).  
Applying our comparison result presented in \cite{38}, it can be proved the existence
of a unique viscosity solution for the equation (\ref{EQ}), with an assigned continuous and bounded initial datum $u_0$, valued 
in the interval $(a,b)$. Moreover, using the same arguments of Theorem 4.5 in \cite{38}, it can be also deduced that, if
the coefficients $\mu$, $\l$, $\eta$, $f$ and the initial datum are $t$-uniformly Lipschitz continuous functions then,
also the viscosity solution has the same regularity at a fixed time. Keeping in mind these features,
we start our considerations from the assumption of the existence of a $t$-uniformly Lipschitz continuous
solution $u$ for the equation (\ref{EQ}) and, then, we prove the semiconvexity/concavity property of it.\\
The first result shows that $u(t)$ is a semiconvex function over the whole space. In the recent literature of viscosity solutions for Hamilton-Jacobi equations, this kind of property
is not still proved for a class of problems which could include our differential model. In \cite{REG}, the authors only prove the semiconcavity of the solution for the Bellman equations,
but the type of nonlinearity of their equation, does not include our Hamiltonian (\ref{H}). 
There is another important paper, by Y. Giga, S. Goto, H. Ishii, and M.H. Sato, \cite{3}, where it was proved
the convexity preserving property of the solution. But also in that case the authors study
a fully nonlinear equation, whose Hamiltonian does not depend on the unknown $u$. 
Therefore, although the initial datum in the original financial model (\ref{DM1}) is constant, their technique
is not suitable for treating our equations (\ref{EQ}), (\ref{DM1}).\\
Using our results about the semiconvexity and semiconcavity and assuming the
same regularity of the coefficient used in (\ref{REG}), we deduce the $\W^{2,\infty}$ regularity, uniformly
in $t\in [0,T)$.\\
In Theorem \ref{SPACETIME}, we conclude next, with a global result of regularity
which contains also the time regularity of the solution. In particular, this last result will follow through an application of the Lipschitz continuity of $u$
with respect to the time variable, which we will present in section 4.\\ 
The general framework for deriving a financial princing equation consists in the application of the classical Ito's rule. Therefore,
we should consider classical solutions which have continuous derivatives. However, there are 
some works about some generalizations of that formula, see for istance in \cite{ITO}. Here it is proved that the Ito's formula holds in arbitrary dimensions for $f\in\W^{2,1,\infty}(\R^N\times(0,T))$, 
if the equation
\begin{eqnarray}\label{ITORULE}
df(X_t,t)={\cal L}f(X_t,t)dt+\nabla^{\top} f(X_t,t)\s(t)\cdot dW_t,\;\;\mbox{a.s.},
\end{eqnarray}
where
\begin{eqnarray*}
{\cal L}f(x,t)&=&\partial_t f(x,t)+\la\mu(x,T-t),\nabla f(x,t)\ra\nonumber\\
&&+\frac{1}{2}tr\big(\s\s^{\top}(T-t)\nabla^2f(x,T-t)\big),\nonumber
\end{eqnarray*}
are interpreted appropriately using the generalized Hessian.\\
In section \ref{Conclusions} we derive some mathematical and financial consequences concerning the model (\ref{DM1}). In
particular we present a result which combines the regularity of the solution of the model (\ref{DM1})
and the results of \cite{ITO} to state the relation (\ref{GABAIX}).
To reach this objective we propose a technical condition about the process
$X_t$ which allows us to extend the Ito's rule as a conditional expected value representation, see Lemma \ref{HAUSSMANN0} later. The technical results proposed in section \ref{Conclusions} allow to consider the degeneration of the process $X_t$, produced by the inequality $N>d$. That degeneration was studied in literature to give sufficient conditions to state the existence of densities for the solutions of stochastic differential equations, see for instance \cite{BELL}, which contain a probabilistic form of H$\ddot o$rmander's. Our condition is of a different nature and applies to many practical cases, like for constant coefficients, where the resultsa of \cite{BELL} do not hold.

\section{Main Results}\label{REG}
In this section we present our main results about the regularity of the solutions, 
the Section 3 being devoted to their proofs.\\
Here we leave out the definition of a viscosity sub/super solution of a differential problem like
(\ref{EQ}), and refer the reader to some classical works about the viscosity theory,
such as \cite{1} or \cite{KR}.\\

\begin{Def}\label{SEMI-CONVEX}
A function $g$ in $\R^N$ is said to be semiconvex with constant $L>0$, if
\begin{eqnarray}\label{SEMI-CONVEX1}
g(x+h)+g(x-h)-2g(x)\geq -L|h|^2,
\end{eqnarray}
for every $x,\;h\in\R^N$.
\end{Def}

In a same way we will say that $g$ is semiconcave with constant $L>0$ if $-g$ is semiconvex
with constant $L$.\\
Over $\R^m$, we will denote $\la\cdot,\cdot\ra$ the usual standard scalar product, while,
over the space, $\M_{m,n}(\R)$, of matrices with real coefficients and $m$-rows and $n$-columns, we will consider the usual norm,
\begin{eqnarray*}
\|A\|=\sup_{x\in \R^n} |Ax|,\;\;\;\forall\;A\in\M_{m,n}(\R).
\end{eqnarray*}
Moreover we will denote as $Im(A)$ the range of the linear map defined by the matrix $A$.\\
If $g:\R^N\times[0,T)\times I\rightarrow\R$, where $I$ is a closed interval in $\R$, then we will use the following notations,
\begin{eqnarray*}
\|g\|_{\infty}&=& esssup_{(x,t,u)\;\in\;\R^N\times [0,T)\times I}|g(x,t,u)|,\\\\
\|g(t)\|_{\infty}&=& esssup_{(x,u)\;\in\;\R^N\times I} |g(x,t,u)|,\;\;\;\forall\;t\in[0,T).
\end{eqnarray*}
We will denote as $\W^{k,\infty}(\R^N)$, $k=1,2$, the usual Sobolev Space of bounded functions with weakly bounded derivatives in $\R^N$ of order less o equal to $k$. 
Moreover $\W^{2,1,\infty}(\R^N\times(0,T))$ denotes the Sobolev space of functions $u$, with weakly
derivatives $\partial_t u,\;\partial_i u,\partial^2_{ij} u\in L^{\infty}(\R^N\times(0,T))$, for $i,j=1,\ldots,N$.
These spaces are respectively endowed with the following norms
\begin{eqnarray*}
\|u\|_{\W^{k,\infty}}&=&\sum_{i=0}^k \|\nabla^i u\|_{\infty},\\\\
\|u\|_{\W^{2,1,\infty}}&=&\|\partial_t u\|_{\infty}+\|\nabla u\|_{\infty}+\|\nabla^2 u\|_{\infty}.
\end{eqnarray*}
If the functions take values in $\R^M$, we refer the same
properties to the single components, and a similar notation holds for functions which depend on another variable in
$I$.\\ 

In the following we shall assume the existence of a viscosity solution $u$ of (\ref{EQ}),
such that $u(\cdot,t)\in\W^{1,\infty}(\R^N)$ and whose norm is uniformly bounded in $t\in[0,T)$.\\

\begin{teo}\label{SEMICONVEX}
Let the viscosity solution $u$ of the equation (\ref{EQ}) be valued in a bound\-ed closed subinterval 
$I$ of the domain $(a,b)$, and let $\mu(\cdot,t),\;w(\cdot,t)\in \W^{2,\infty}(\R^N)$, 
$f(\cdot,t,\cdot)\in\W^{2,\infty}(\R^N\times I)$, uniformly in time.\\
Assume:
\begin{description}
\item[$i$)] $\l\in\C((a,b))$, $\eta\in\C((a,b))\cap\C^2(I)$;
\item[$ii$)] For every $(x,t)\in\R^N\times[0,T)$, $w(x,t)\in Im(\s^{\top}(t))$.
\end{description}
If $u_0$ is semiconvex, then there are positive constants $C,\;M_0,\;C_0$, such that 
\begin{eqnarray}\label{CHG0}
u(x+h,t)+u(x-h,t)-2 u(x,t)\geq -(M_0 e^{C t}+C_0)|h|^2,
\end{eqnarray}
holds for every $x,\;h\in\R^N$, and $t\in[0,T)$. Therefore, for every $t\in[0,T)$, $u(\cdot,t)$ is 
semiconvex.
\end{teo}

We have also the analogous result for the semiconcavity of $u$.

\begin{teo}\label{SEMICONCAVE}
Let the viscosity solution $u$ of the equation (\ref{EQ}) be valued in a bound\-ed closed subinterval 
$I$ of the domain $(a,b)$, and let $\mu(\cdot,t),\;w(\cdot,t)\in\W^{2,\infty}(\R^N)$, 
$f(\cdot,t,\cdot)\in\W^{2,\infty}(\R^N\times I)$, uniformly in time.
Assume:
\begin{description}
\item[$i$)] $\l\in\C((a,b))$, $\eta\in\C((a,b))\cap\C^2(I)$;
\item[$ii$)] For every $(x,t)\in\R^N\times[0,T)$, $w(x,t)\in Im(\s^{\top}(t))$.
\end{description}
If $u_0$ is semiconcave, then there are positive constants $C,\;M_0,\;C_0$, such that 
\begin{eqnarray}\label{CHG0}
u(x+h,t)+u(x-h,t)-2 u(x,t)\leq (M_0 e^{C t}+C_0)|h|^2,
\end{eqnarray}
holds for every $x,\;h\in\R^N$, and $t\in[0,T)$. Therefore, for every $t\in[0,T)$, $u(\cdot,t)$ is 
semiconcave.
\end{teo}

As a consequence, we will deduce the announced result about the regularity of the solution $u$.

\begin{teo}\label{regf1}
If $\mu$, $w$, $f$, $\l$ and $\eta$ satisfy the assumptions of Theorems \ref{SEMICONVEX} and \ref{SEMICONCAVE} and
$u_0\in\W^{2,\infty}(\R^N)$, then $u(t)\in\W^{2,\infty}(\R^N)$, uniformly in time. 
\end{teo}

\begin{teo}\label{SPACETIME}
If $\s$ is Lipschitz continuous, $\mu,\;w\in\W^{2,1,\infty}(\R^N\times(0,T))$, $f\in\W^{2,1,2,\infty}(\R^N\times(0,T)\times I)$,
$w(x,t)\in Im(\s^{\top}(t))$ for every $(x,t)\in\R^N\times[0,T)$, $u_0\in\W^{2,\infty}(\R^N)$, then, 
\begin{eqnarray}\label{SPACETIME2}
u\in\W^{2,1,\infty}(\R^N\times(0,T)).
\end{eqnarray}
\end{teo}

To prove Theorems \ref{SEMICONVEX}, \ref{SEMICONCAVE} and \ref{regf1}, we need of some
technical results, see the Propositions \ref{regularity} and \ref{regn}, in section 3, which
state the same thesis but use some structural and regularity assumptions on the function $\l$.\\
Then, these conditions on $\l$, can be removed by a particular compatibility between the second order linear term and the first
order term in the equation (\ref{EQ}). The proof of these facts is based on the classical property about the preservation of
the notion of viscosity solution through a global, increasing, change of the variable $u$.\\

\section{Proof of the Results}\label{PROOFS}

This part is devoted to the presentation of the technical results which are useful for proving Theorems \ref{SEMICONVEX}, \ref{SEMICONCAVE} and \ref{regf1}, illustrated in the previous section.

\begin{prop}\label{regularity}
Consider a viscosity solution $u$ of problem (\ref{EQ}) valued in a bound\-ed closed subinterval $I$ of the domain $(a,b)$, such that $u(\cdot,t)\in\W^{1,\infty}(\R^N)$, with a norm uniformly bounded in the time $t\in[0,T)$. 
Assume a semiconvex initial datum $u_0$ (see definition \ref{SEMI-CONVEX}) with a constant $L_0>0$.
Suppose that, $\mu(\cdot,t),\;w(\cdot,t)$ are $\W^{2,\infty}(\R^N)$ functions and, $f(\cdot,t,\cdot)\in\W^{2,\infty}(\R^N\times I)$, uniformly in time. Moreover 
assume:
\begin{description}
\item[$i$)] $\l,\;\eta\in\C((a,b))\cap\C^2(I)$. \\
\item[$ii$)] $\l<0$, $\l^{\p}>0$, $\l\l^{\p\p}-2(\l^{\p})^2>0$, over $I$.\\
\item[$iii$)] For every $(x,t)\in\R^N\times[0,T)$, $w(x,t)\in Im(\s^{\top}(t))$.
\end{description}
Then, there are positive constants $C$, which depends on $\sup_{t\in[0,T)}\|u(t)\|_{\W^{1,\infty}}$
and $C_0$, which depends on $L_0$ and $Lip(u_0)$, such that
\begin{eqnarray}\label{EQassertion}
u(x+h,t)+u(x-h,t)-2 u(x,t)\geq -C_0 e^{C t} |h|^2
\end{eqnarray}
holds for every $(x,h,t)\in\R^N\times\R^N\times[0,T)$. In particular $u(\cdot,t)$ is semiconvex, for every $t\in[0,T)$.
\end{prop}

The equivalent result for the semiconcave property of the solution, is the following.

\begin{prop}\label{regn}
Consider a viscosity solution $u$ of problem (\ref{EQ}) valued in a bound\-ed closed subinterval, $I$ of the domain $(a,b)$, such that $u(\cdot,t)\in\W^{1,\infty}(\R^N)$, with a norm uniformly bounded in the time $t\in[0,T)$. 
Assume a semiconcave initial datum $u_0$ (see definition \ref{SEMI-CONVEX}) with a constant $L_0>0$.
Suppose that, $\mu(\cdot,t),\;w(\cdot,t)$ are $\W^{2,\infty}(\R^N)$ functions and, $f(\cdot,t,\cdot)\in\W^{2,\infty}(\R^N\times I)$, uniformly in time. Moreover 
assume:
\begin{description}
\item[$i$)] $\l,\;\eta\in\C((a,b))\cap\C^2(I)$. \\
\item[$ii$)] $\l>0$, $\l^{\p}>0$, $\l\l^{\p\p}-2(\l^{\p})^2>0$, over $I$.\\
\item[$iii$)] For every $(x,t)\in\R^N\times[0,T)$, $w(x,t)\in Im(\s^{\top}(t))$.
\end{description} 
Then, there are positive constants $C$, which depends on $\sup_{t\in[0,T)}\|u(t)\|_{\W^{1,\infty}}$
and $C_0$, which depends on $L_0$ and $Lip(u_0)$, such that
\begin{eqnarray}\label{EQ1}
u(x+h,t)+u(x-h,t)-2 u(x,t)\leq C_0 e^{C t} |h|^2
\end{eqnarray}
holds for every $(x,h,t)\in\R^N\times\R^N\times[0,T)$. In particular $u(\cdot,t)$ is semiconcave, for every $t\in[0,T)$.
\end{prop}

\begin{oss}\rm
We can remark the analogies between the two formulations \ref{regularity}, \ref{regn}. Actually the hypothesis
$\l^{\p}>0$, is the same in both the Propositions. This can be interpreted as a consequence of the comparison principle,
which requires a monotonicity for the hamiltonian with respect to $u$, (see \cite{3}, pag. 462).
While of course, the others expressions in $ii$) have exactly an opposite sign, showing a perfect 
symmetry between the two formulations.
\end{oss}

The proof of Proposition \ref{regn} follows by argues as Proposition \ref{regularity}
and then it is omitted.\\
The proof is based on the same technique proposed in \cite{REG}, which uses
a particular test function. Nevertheless the nonlinear part in the equation (\ref{EQ}) does not 
allow to close the proof as in that work. Hence, we turn to a semiconcavity property of the nonlinear part with respect to $(u,\nabla u)$, and to
a monotonicity with respect to the unknown.\\  

\begin{oss}\label{U0}\rm
We observe that if a function $g$ is a semiconvex as in the Definition \ref{SEMI-CONVEX} and is a Lipschitz continuous function, then, for all $x,y,z$, the following inequality holds, 
\begin{eqnarray*}
g(x)+g(y)-2g(z)\geq -L(|x-z|^2+|y-z|^2)-Lip(g)|x+y-2z|.
\end{eqnarray*} 
Where $Lip(g)$ is the Lipschitz constant of $g$.
 \end{oss}

\begin{Dimo} {\bf Proposition \ref{regularity}.}
Consider the function $v=u e^{-Ct}$, where $C$ is a nonnegative constant.
Then $v$ is an $x$-Lipschitz continuous function uniformly with respect to the time $t\in [0,T)$. Moreover $v$ is 
a continuous viscosity solution of the following equation

\begin{eqnarray}\label{EQ3}
\partial_t v-\frac{1}{2}tr(\s\s^{\top} \nabla^2 v)+\la \mu, \nabla v \ra+\l(v e^{C t})e^{Ct}|\s^{\top} \nabla v|^2\nonumber\\
+\eta(v e^{Ct})\la \s^{\top} v, w \ra +e^{-Ct}f(x,t,v e^{Ct})+ C v=0,
\end{eqnarray}


where $(x,t)\in \R^N\times (0,T)$, and $v_0=u_0$. We now will go to prove that, under the assumptions made on 
$\l$ and $\eta$, we have the semiconvexity of the function $v$, and therefore the semiconvexity
of the function $u$. In particular we are going to show that, 

\begin{eqnarray}\label{EQinequality}
v(x,t)+v(y,t)-2v(z,t)\geq -M(|x-z|^4+|y-z|^4+|x+y-2z|^2)^{\frac{1}{2}}
\end{eqnarray}

for every $x,y,z\in \R^N,\;t\in [0,T)$, where $M=\frac{\sqrt{3}}{2}\max(L_0,Lip(u_0))$. By the assumptions on $u$, the initial datum is Lipschitz continuous, therefore this constant is well defined; this obviously yields the assertion on $v$, by pluging $x=z+h$ and $y=z-h$. It is easy to see that the above inequality is equivalent to the following
one:

\begin{eqnarray}\label{EQequiv}
\lefteqn{v(x,t)+v(y,t)-2v(z,t)\geq}\nonumber\\ 
&&- M \big[\delta+\frac{1}{\delta}(|x-z|^4+|y-z|^4+|x+y-2z|^2)\big],
\end{eqnarray}
for every $\delta>0,\;\forall\;x,y,z \in\R^N $.


Hence, fix $\var,\delta,\gamma>0$, and consider the following test function


\begin{eqnarray}\label{EQtest}
&&\Psi(x,y,z,t)= v(x,t)+v(y,t)-2v(z,t)+M \big[\delta+\frac{1}{\delta}(|x-z|^4\nonumber\\ &&+|y-z|^4+|x+y-2z|^2)\big]
+\var|x|^2+\var|y|^2+\var|y|^2+\frac{\gamma}{T-t},
\end{eqnarray}


defined for $(x,y,z,t)\in U=\R^{3N}\times [0,T)$. The assertion (\ref{EQequiv}) is equivalent to proving
that, for every $\delta,\gamma>0$, there exists $\var_0=\var_0(\delta,\gamma)>0$, such that
for every $0<\var<\var_0$, the following holds:  


\begin{eqnarray}\label{EQequiv1}
\inf_{U} \Psi\geq 0.
\end{eqnarray}


Actually if (\ref{EQequiv1}) holds, then fixing a point over $U$, we can send $\var$ to zero in the inequality $\Psi\geq 0$, and then
sending also $\gamma$ to zero we obtain (\ref{EQequiv}). Thus we limit us to consider the assertion (\ref{EQequiv1}).\\
We assume as usual that (\ref{EQequiv1}) is false, and will get a contradiction. Therefore
exist $\delta_0,\;\gamma_0>0$ and a sequence $\var_j\rightarrow 0$, as $j\rightarrow \infty$, such
that 


\begin{eqnarray}\label{EQequiv2}
\inf_{U} \Psi<0,
\end{eqnarray}


with $\delta=\delta_0$, $\gamma=\gamma_0$, and $\var=\var_j$, for every integer $j>0$. If we consider a minimizing sequence $(x_k,y_k,z_k,t_k)\in U$ for $\Psi$. By (\ref{EQequiv2}), the definition (\ref{EQtest}) and the boundness of $v$, we see that $(x_k,y_k,z_k)$ must be bounded, so we can extract a convergent subsequence, which converges to some point $(\w{x},\w{y},\w{z},\w{t})\in U$, which, by the continuity of $v$, is a global minimum point for $\Psi$ over $U$; moreover $\w{t}$ (which is obviously less than $T$) is strictly

positive. In fact if $\w{t}=0$, by (\ref{EQequiv2}) and the Remark \ref{U0}, we obtain,


\begin{eqnarray*}
0&>&\Psi(\w{x},\w{y},\w{z},0) \geq \nonumber\\
&\geq&u_0(\w{x})+u_0(\w{y})-2u_0(\w{z})+M \big[\delta+\frac{1}{\delta}(|\w{x}-\w{z}|^4+|\w{y}-\w{z}|^4+|\w{x}+\w{y}-2\w{z}|^2)\big]\geq\\ 
&\geq&-L_0 (|\w{x}-\w{z}|^2+|\w{y}-\w{z}|^2)-Lip(u_0)|\w{x}+\w{y}-2\w{z}|+\\
&&+M\big[\delta+\frac{1}{\delta}(|\w{x}-\w{z}|^4+|\w{y}-\w{z}|^4+|\w{x}+\w{y}-2\w{z}|^2)\big]\geq\\
&\geq& [-L_0+\frac{2}{\sqrt{3}}M] (|\w{x}-\w{z}|^2+|\w{y}-\w{z}|^2)\nonumber\\
&&+ [-Lip(u_0)+\frac{2}{\sqrt{3}}M]|\w{x}+\w{y}-2\w{z}|\geq 0.
\end{eqnarray*}
where in the last two inequalities we have used the following relations

\begin{eqnarray*}
\sqrt{x_1^2+x_2^2+x_3^2}&\geq& \frac{1}{\sqrt{3}}(x_1+x_2+x_3),\;\;\;\forall\;x_1,x_2,\;x_3\geq 0,\\
x_1^2+x_2^2&\geq& 2 x_1 x_2,\;\;\;\forall x_1,x_2 \in \R
\end{eqnarray*}


Hence by the previous contradiction we deduce that $\w{t}>0$. So the minimum point is an interior stationary point
of $\Psi$. Setting the function

\begin{eqnarray}\label{EQfi}
g(x,y,z,t)&=&-v(x,t)-v(x,t)+2 v(z,t),\\
\Phi(x,y,z,t)&=& M \big[\delta+\frac{1}{\delta}(|x-z|^4+|y-z|^4+|x+y-2z|^2)\big]\nonumber\\
&&+\var|x|^2+\var|y|^2+\var|y|^2+\frac{\gamma}{T-t}
\end{eqnarray}


we have that $g-\Phi=-\Psi$ has a global interior maximum point at $(\w{\xi},\w{t})=(\w{x},\w{y},\w{z},\w{t})$, therefore we can apply the classical
Theorem about the maximum principle for semicontinuous functions of M.G. Crandall and H. Ishii,   
in \cite{2}, to deduce that for $\kappa=\frac{1}{\var}>0$ there exist $(b_i,X_i)\in\R\times \S^N$, for $i=1,2,3$, such
that,


\begin{eqnarray}\label{EQappartz}
&&(-b_i,-\Phi_i, -X_1)\in \o{P}^{2,-}v_i,\;\;i=1,2\\
&&(\frac{b_3}{2}, \frac{\Phi_3}{2}, \frac{1}{2}X_3)\in \o{P}^{2,+}v_3,
\end{eqnarray}


and, if $O$ denotes the null $N\times N$ matrix, we have 


\begin{eqnarray}\label{EQappartz1}
\left(\begin{array}{ccc}
X_1 & 0   & 0\\
0   & X_2 & 0\\
0   & 0   & X_3
\end{array}
\right)&\leq& \nabla^2 \Phi(\w{\xi},\w{t})+\kappa [ \nabla^2 \Phi(\w{\xi},\w{t})]^2,\\
b_1+b_2+b_3 &=& \partial_t \Phi(\w{\xi},\w{t}).
\end{eqnarray}


Where for simplifying notations we have set $v_1,v_2,v_3$ for $v(\w{x},\w{t})$, $v(\w{y},\w{t})$, $v(\w{z},\w{t})$, respectively,
and in a same way $\Phi_1,\Phi_2,\Phi_3$, for the partial derivatives of $\Phi$ evaluated at the considered maximum point.
Now we compute the derivatives of $\Phi$. Set 


\begin{eqnarray}\label{EQab}
p &=&\frac{2 M}{\delta}(\w{x}-\w{z})|\w{x}-\w{z}|^2\\
q &=&\frac{2 M}{\delta}(\w{y}-\w{z})|\w{y}-\w{z}|^2\\
m &=&\frac{2 M}{\delta}(\w{x}+\w{y}-2\w{z}).
\end{eqnarray}
%

Then


\begin{eqnarray}\label{EQderivatives}
\partial_t \Phi(\w{\xi})&=&\frac{\gamma}{(T-\w{t})^2}\geq \frac{\gamma}{T^2}\nonumber\\ \nonumber\\
\Phi_1 &=& 2\var \w{x}+2 p+m \nonumber\\ 
\Phi_2 &=& 2\var \w{y}+2 q+m  \nonumber\\
\Phi_3  &=& 2\var \w{z}-2p-2q-2m.
\end{eqnarray}
%

Moreover, if $I$ denotes the $N\times N$ identity matrix, we have 


\begin{eqnarray}\label{EQmad}
\nabla^2 \Phi(\w{\xi},\w{t})&=&
(\frac{2M}{\delta})^{\frac{1}{3}}\Big\{ 
2\left(\begin{array}{ccc} |p|^{\frac{2}{3}}I & 0 & -|p|^{\frac{2}{3}}I \\
   0 & |q|^{\frac{2}{3}}I & -|q|^{\frac{2}{3}}I \\
           -|p|^{\frac{2}{3}}I & -|q|^{\frac{2}{3}}I & (|p|^{\frac{2}{3}}+|q|^{\frac{2}{3}})I
\end{array}\right)\nonumber\\
&&+ 4 \left(\begin{array}{ccc}
\frac{p\otimes p}{|p|^{\frac{4}{3}}} & 0 & -\frac{p\otimes p}{|p|^{\frac{4}{3}}}\\
0 &  \frac{q\otimes q}{|q|^{\frac{4}{3}}} & - \frac{q\otimes q}{|q|^{\frac{4}{3}}}\\
- \frac{p\otimes p}{|p|^{\frac{4}{3}}} & - \frac{q\otimes q}{|q|^{\frac{4}{3}}} &  \frac{p\otimes p}{|p|^{\frac{4}{3}}}+\frac{q\otimes q}{|q|^{\frac{4}{3}}}
\end{array}\right)\Big\}\nonumber\\
&&+\frac{2 M}{\delta}
\left(\begin{array}{ccc}
I & I & -2 I \\
I & I & -2 I \\
-2I & -2I  & 4I 
\end{array}\right)+ 2\var 
\left(\begin{array}{ccc}
I & 0 & 0 \\
0 & I & 0 \\
0 & 0 & I 
\end{array}\right)
\end{eqnarray}


By (\ref{EQappartz}), and the equation (\ref{EQ3}), we have


\begin{eqnarray}\label{EQsubsti}
b_i &\leq& \frac{1}{2}tr(\s\s^{\top}X_i)-\la \mu_i, \Phi_i \ra +\l(v_i e^{C \wt})e^{C \wt }|\s^{\top}\Phi_i|^2\nonumber\\
&&-\eta (v_i e^{C \wt})\la \s^{\top} \Phi_i, w_i \ra+C v_i +e^{-C \wt}f_i,\;\;\;\;\;\;\;\;\;\;\;\;i=1,2 \nonumber\\
b_3 &\leq &\frac{1}{2}tr(\s\s^{\top}X_3)-\la \mu_3, \Phi_3 \ra-2\l(v_3 e^{C \wt})e^{C \wt}|\frac{1}{2}\s^{\top}\Phi_3|^2\nonumber\\
&&-\eta(v_3 e^{C \wt})\la \s^{\top}\Phi_3,w_3\ra -2e^{-C \wt}f_3-2C v_3.
\end{eqnarray}
%

Where we have used the same notation of (\ref{EQappartz}) for denoting the functions $\mu,w,f$, and we have also omitted 
the dipendency of the matrix $\s$ by $\w{t}$. Adding inequalities (\ref{EQsubsti}) and using
(\ref{EQappartz1}), we compute


\begin{eqnarray}\label{EQadding}
\partial_t\Phi(\w{\xi},\wt)&\leq&
\frac{1}{2}tr\big(\s\s^{\top}(X_1+X_2+X_3)\big)-\big[
\la \mu_1, \Phi_1 \ra+\la \mu_2, \Phi_2 \ra+\la \mu_3, \Phi_3 \ra\big]+\nonumber\\
&&+\big[\l(v_1 e^{C \wt})e^{C \wt }|\s^{\top}\Phi_1|^2
-\eta (v_1 e^{C \wt})\la \s^{\top} \Phi_1, w_1 \ra\nonumber\\
        &&+\l(v_2 e^{C \wt})e^{C \wt }|\s^{\top}\Phi_2|^2
-\eta (v_2 e^{C \wt})\la \s^{\top} \Phi_2, w_2 \ra\nonumber\\
&&-2\l(v_3 e^{C \wt})e^{C \wt}|\frac{1}{2}\s^{\top}\Phi_3|^2
-\eta(v_3 e^{C \wt})\la \s^{\top}\Phi_3,w_3\ra\big]+\nonumber\\
&&+C(v_1+v_2-2v_3)+ e^{-C\wt}[f_1+f_2-2 f_3].
\end{eqnarray} 


Now we use (\ref{EQderivatives}) and (\ref{EQmad}) for estimating the single part in the brackets of the inequality
(\ref{EQadding}). Let $\s^{(l)}$ be the $l$ column of $\s$, and define


\[\Sigma^{(l)}=\left(\begin{array}{c}

\s^{(l)} \\

\s^{(l)} \\

\s^{(l)}

\end{array}\right)\]


for $l=1,\ldots,d$, then by (\ref{EQappartz1}), we have


\begin{eqnarray}\label{EQlinear2}
tr\big(\s\s^{\top}(X_1+X_2+X_3)\big)&=&\sum_{l=1}^d\big(\la X_1\s^{(l)}, \s^{(l)} \ra
+\la X_2\s^{(l)}, \s^{(l)} \ra+\la X_3\s^{(l)}, \s^{(l)} \ra\big)\nonumber\\
&\leq& \sum_{l=1}^d\big(\la \nabla^2 \Phi(\w{\xi},\wt) \Sigma^{(l)}, \Sigma^{(l)}\ra+\nonumber\\
&& \kappa |\nabla^2 \Phi(\w{\xi},\wt) \Sigma^{(l)}|^2\big).
\end{eqnarray}


Using (\ref{EQmad}), it is easy to compute that 


\begin{eqnarray}\label{EQlinear21}
\nabla^2\Phi(\w{\xi},\wt) \Sigma^{(l)}=2 \var \Sigma^{(l)}.
\end{eqnarray}


So introducing (\ref{EQlinear21}) in (\ref{EQlinear2}), we deduce


\begin{eqnarray}\label{EQlinear22}
tr\big(\s\s^{\top}(X_1+X_2+X_3)\big)\leq 18 \var\sum_{l=1}^d |\s^{(l)}|^2=
O(\var).
\end{eqnarray}


By condition (\ref{EQequiv2}), we have 


\begin{eqnarray}\label{EQxy}
&&\var |\wx|, \var |\wy|, \var |\wz| \leq  O(\sqrt{\var}),\nonumber\\
\\
&& \frac{1}{2}[(\frac{\delta}{2 M})^{\frac{1}{3}}(|p|^{\frac{4}{3}}+|q|^{\frac{4}{3}})+
\frac{\delta}{2 M} |m|^2]+M\delta\leq 2 v_3-v_1-v_2.\nonumber
\end{eqnarray} 


Therefore, using the notations introduced in (\ref{EQab})-(\ref{EQderivatives}), (\ref{EQxy}) and the regularity assumptions on $\mu$, we obtain


\begin{eqnarray}\label{EQlinear1}
\lefteqn{|\la \mu_1, \Phi_1 \ra + \la \mu_2, \Phi_2 \ra + \la \mu_3, \Phi_3 \ra |\leq}\nonumber\\
&\leq&|\la \mu_1, 2 p+m \ra + \la  \mu_2, 2 q+m  \ra - \la \mu_3, 2p+2q+2m \ra |+\nonumber\\
&&+2\var |\la \mu_1, \wx \ra + \la \mu_2, \wy \ra + \la \mu_3, \wz \ra |\leq\nonumber\\
&\leq&| 2 \la\mu_1 - \mu_3, p \ra +2 \la \mu_2 - \mu_3, q \ra + \la \mu_1 + \mu_2 - 2 \mu_3, m\ra|
+O(\sqrt{\var})\leq \nonumber\\
&\leq&\sup_{t\in [0,T)}\|\mu (t) \|_{\W^{2,\infty}}\big[2(\frac{\delta}{2 M})^{\frac{1}{3}}(|p|^{\frac{4}{3}}+
|q|^{\frac{4}{3}})+\nonumber\\
&&+[(\frac{\delta}{2 M})^{\frac{1}{3}}(|p|^{\frac{4}{3}}+|q|^{\frac{4}{3}})+
\frac{\delta}{2 M} |m|^2]^{\frac{1}{2}}(\sqrt{\frac{\delta}{2M}}|m|)\big]+O(\sqrt{\var})
\leq \nonumber\\
&\leq& \sup_{t\in [0,T)}\|\mu(t)\|_{\W^{2,\infty}}\Big[
4(2 v_3-v_1-v_2)+\frac{1}{2}
[(\frac{\delta}{2 M})^{\frac{1}{3}}(|p|^{\frac{4}{3}}+|q|^{\frac{4}{3}})+
\frac{\delta}{2 M} |m|^2]+\nonumber\\
&&+\frac{\delta}{4 M}|m|^2\Big]+O(\sqrt{\var})
\leq\nonumber\\
&\leq& 6 \sup_{t\in [0,T)}\|\mu(t)\|_{\W^{2,\infty}}(2 v_3-v_1-v_2)+O(\sqrt{\var}), 
\end{eqnarray}


where in the last inequality we have again used (\ref{EQxy}).








For estimating the nonlinear part we consider the function


\begin{eqnarray}\label{EQnonlinear}
G(u,\theta)=\l(u)e^{C \wt} |\s^{\top}\theta|^2-\eta(u)\la \s^{\top} \theta , w_3\ra, 
\end{eqnarray}


which depends on $(u,\theta)\in [a,b]\times \R^N$, and for every $u_1,u_2,\theta_1,\theta_2$, define


\begin{eqnarray}\label{EQPP}
\Delta^2 G(u_1,u_2,\theta_1,\theta_2)=G(u_1,\theta_1)+G(u_2,\theta_2)-2G(\frac{u_1+u_2}{2},\frac{\theta_1+\theta_2}{2}).
\end{eqnarray}
%

Let $\eta_i$ denotes the value of $\eta$ at the point $v_i e^{C\wt}$, for $i=1,2,3$, then using notations (\ref{EQnonlinear}), (\ref{EQPP}) and (\ref{EQab})-(\ref{EQderivatives}), we have

%
\begin{eqnarray}\label{EQnonlinear1}
\lefteqn{\Big[\l(v_1 e^{C \wt})e^{C \wt }|\s^{\top}\Phi_1|^2
-\eta (v_1 e^{C \wt})\la \s^{\top} \Phi_1, w_1 \ra+\l(v_2 e^{C \wt})e^{C \wt }|\s^{\top}\Phi_2|^2}\nonumber\\
&&-\eta (v_2 e^{C \wt})\la \s^{\top} \Phi_2, w_2 \ra
-2\l(v_3 e^{C \wt})e^{C \wt}|\frac{1}{2}\s^{\top}\Phi_3|^2
-\eta(v_3 e^{C \wt})\la \s^{\top}\Phi_3,w_3\ra\big]=\nonumber\\
&=&\Delta^2 G(v_1 e^{C \wt},v_2 e^{C \wt},2p+m,2q+m)+2\big[G(\frac{(v_1+v_2)e^{C \wt}}{2}, p+q+m)-\nonumber\\
&&-G(v_3 e^{C \wt}, p+q+m)\big]
+\big[\eta_1\la \s^{\top}(2p+m), w_3-w_1\ra+\nonumber\\
&&+\eta_2\la \s^{\top}(2q+m), w_3-w_2\ra\big]
+O(\sqrt{\var}).
\end{eqnarray}
%
In the last passages we have used the inequalities (\ref{EQxy}) for
estimating the residual terms which involve $\var$.

\begin{eqnarray}\label{EQnonlinear2}
\lefteqn{\eta_1\la \s^{\top}(2p+m), w_3-w_1\ra+\eta_2\la \s^{\top}(2q+m), w_3-w_2\ra=}\nonumber\\
&&=[\eta_1-\eta_3]\la \s^{\top}(2p+m), w_3-w_1\ra+[\eta_2-\eta_3]\la \s^{\top}(2q+m), w_3-w_2\ra\nonumber\\
&&+2\eta_3\big[\la \s^{\top}p, w_3-w_1\ra+\la \s^{\top}q, w_3-w_2\ra\big]\nonumber\\
&&-\eta_3\la \s^{\top}m,w_1+w_2-2 w_3\ra.
\end{eqnarray} 


Hence, using the Lipschitz regularity of the function $u$, the regularity of $w$, jointly with (\ref{EQxy}), we deduce,


\begin{eqnarray}\label{EQnonlinear3}
\lefteqn{[\eta_1-\eta_3]\la \s^{\top}(2p+m), w_3-w_1\ra\leq}\\ 
&\leq&\|\eta^{\p}\|_{\infty}\|\s^{\top}\|_{\infty}\sup_{t\in[0,T)}\|w(t)\|_{\W^{2,\infty}}
\big[4\|u\|_{\infty}(\frac{\delta}{2M})^{\frac{1}{3}}|p|^{\frac{4}{3}}+\nonumber\\
&&+\sup_{t\in[0,T)}\|u(\cdot,t)\|_{\W^{1,\infty}} 
(\frac{\delta}{2M})^{\frac{2}{3}}|p|^{\frac{2}{3}}|m|\big]
\leq \nonumber\\
&\leq&\|\eta^{\p}\|_{\infty}\|\s^{\top}\|_{\infty}\sup_{t\in[0,T)}\|w(t)\|_{\W^{2,\infty}}
\big[4\|u\|_{\infty}(\frac{\delta}{2M})^{\frac{1}{3}}|p|^{\frac{4}{3}}+
\frac{1}{2} \sup_{t\in[0,T)}\|u(t)\|_{\W^{1,\infty}}\nonumber\\ 
&&((\frac{\delta}{2M})^{\frac{1}{3}}|p|^{\frac{4}{3}}+\frac{\delta}{2 M}|m|^2)\big]
\leq C_1\big[\frac{1}{2} (\frac{\delta}{2M})^{\frac{1}{3}}|p|^{\frac{4}{3}}+\frac{1}{4}
\frac{\delta}{2 M}|m|^2)\big].
\end{eqnarray}


In the last inequality we have used the notation


\begin{eqnarray}\label{EQconstant1}
C_1=2\|\eta^{\p}\|_{\infty}\|\s^{\top}\|_{\infty}\sup_{t\in[0,T)}\|w(t)\|_{\W^{2,\infty}}
\max(8\|u\|_{\infty}, \sup_{t\in[0,T)}\|u(t)\|_{\W^{1,\infty}})
\end{eqnarray}


We can repeat the argument for estimating $[\eta_2-\eta_3]\la \s^{\top}(2q+m), w_3-w_2\ra$, obtaining, again using the second relation in (\ref{EQxy}), the inequality


\begin{eqnarray}\label{EQnonlinear4}
\lefteqn{[\eta_1-\eta_3]\la \s^{\top}(2p+m), w_3-w_1\ra+[\eta_2-\eta_3]\la \s^{\top}(2q+m), w_3-w_2\ra
\leq}\nonumber\\
&\leq& C_1\big[\frac{1}{2} \frac{\delta}{2M})^{\frac{1}{3}}(|p|^{\frac{4}{3}}+|q|^{\frac{4}{3}})+
\frac{1}{2}\frac{\delta}{2 M}|m|^2)\big]\leq C_1(2 v_3-v_1-v_2)
\end{eqnarray}


By (\ref{EQxy}) and the following inequality,


\begin{eqnarray*}
|w_1+w_2-2w_3|\leq \sup_{t\in[0,T)}\|w(t)\|_{\W^{2,\infty}}\big[\frac{\delta |p|}{2M})^{\frac{4}{3}}+
\frac{\delta |q|}{2M})^{\frac{4}{3}}+\frac{\delta |m|}{2M})^{2}\big]^{\frac{1}{2}},
\end{eqnarray*}


then it is easy to estimate the last two terms in the brackets $[\cdot]$ in (\ref{EQnonlinear2}), through the expression $C_2(2v_3-v_1-v_2)$, where


\begin{eqnarray}\label{EQconstant2}
C_2=5\|\eta\|_{\infty}\|\s^{\top}\|_{\infty}\sup_{t\in[0,T)}\|w(t)\|_{\W^{2,\infty}}.  
\end{eqnarray}


Now we proceed with the estimates of the first two terms in (\ref{EQnonlinear1}). We observe that, by the assumption $ii$), we can consider the positive constant,


\begin{eqnarray}\label{EQconstant3}
C_3=\frac{1}{2}\frac{\|\eta^{\p}\|^2_{\infty}\|w\|^2_{\infty}}{4 \min \l^{\p}},
\end{eqnarray}


and also we can write


\begin{eqnarray}\label{EQnonlinear5}
\partial_u G(u,\theta)\geq - C_3 e^{-C \wt},\;\;\;\forall\;(u,\theta)\in[a,b]\times\R^N.
\end{eqnarray}


Therefore by (\ref{EQnonlinear5}) and $2v_3-v_1-v_2>0$, the following holds,


\begin{eqnarray}\label{EQnonlinear6}
G(\frac{(v_1+v_2)e^{C \wt}}{2}, p+q+m)-G(v_3 e^{C \wt}, p+q+m)\leq C_3(2 v_3-v_1-v_2).
\end{eqnarray}


By $iii$), we can set $w_3=\s^{\top}b_3$, for some $N$-dimensional vector $b_3$. Moreover


\begin{eqnarray}\label{EQnonlinear7}
\partial_u^{2}G &=& e^{C \wt}\l^{\p\p}|\s^{\top}\theta|^2-\eta^{\p\p}\la\s^{\top}\theta,w_3\ra, \nonumber\\
\partial_{\theta,u}^2 G &=& \s\s^{\top}J,\nonumber\\
\partial_{\theta}^2 G &=& 2 e^{C \wt}\l \s\s^{\top},
\end{eqnarray}


where $J=2 e^{C \wt} \l^{\p}\theta-\eta^{\p} b_3$. For every $X=(k,h)\in\R\times\R^N$ consider the orthogonal matrix $A$ such that $A\s\s^{\top}A^{\top}$ is diagonal, with entries $S_i\geq 0$, $i=1,\ldots,N$. Set $\o{h}=A h$, $\o{J}=AJ$, and define


\begin{eqnarray}\label{EQconstant4}
C_4=\frac{\|\eta^{\p\p}-2\l^{\p}\eta^{\p}\|_{\infty}^2 \|w\|_{\infty}^2}{4 \min(2\frac{(\l^{\p})^2}{\l}-\l^{\p\p})}
     -\frac{1}{2\max(\l)}\|\eta^{\p}\|_{\infty}^2\|w\|_{\infty}^2,
\end{eqnarray}


by (\ref{EQnonlinear7}), $ii$), and (\ref{EQconstant4}), we have


\begin{eqnarray}\label{EQnonlinear8}
\lefteqn{\la \left(\nabla^{2} G-C_4 e^{-C\wt}\left(\begin{array}{cc}
I & 0\\
0 & 0\end{array}\right)\right)X,X\ra= }\nonumber\\
&=& (\partial_u^2 G-C_4 e^{-C\wt})k^2+\sum_{i=1}^{N}(2 e^{C \wt}\l \o{h}_i^2+2 k\o{J}_i\o{h}_i)S_i 
\leq\nonumber\\
&\leq& \big[\partial_u^2 G-C_4 e^{-C\wt}-\sum_{i=1}^N \frac{e^{-C \wt}}{2\l}S_i|\o{J}_i|^2\big]k^2
=\nonumber\\&=&\big[
e^{C \wt}\l^{\p\p}|\s^{\top}\theta|^2+\eta^{\p\p}\la\s^{\top}\theta,w_3\ra,
-C_4 e^{-C\wt}-\frac{e^{-C \wt}}{2\l}|2e^{C \wt} \l^{\p}\s^{\top}\theta+\eta^{\p} w_3|^2\big]k^2\nonumber\\
&=&\big[e^{C \wt}(\l^{\p\p}-2\frac{(\l^{\p})^2}{\l})|\s^{\top}\theta|^2+
(\eta^{\p\p}-2\l^{\p}\eta^{\p})\la\s^{\top}\theta,w_3\ra-\nonumber\\ 
&&-\frac{e^{-C \wt}}{2\l}|\eta^{\p}|^2|w_3|^2-C_4 e^{-C\wt}\big]k^2
\leq \nonumber\\
&\leq& \big[e^{-C \wt}\frac{\|\eta^{\p\p}-2\l^{\p}\eta^{\p}\|_{\infty}^2 \|w\|_{\infty}^2}{4 \min(2\frac{(\l^{\p})^2}{\l}-\l^{\p\p})}\nonumber\\
&&-\frac{e^{-C \wt}}{2\max(\l)}\|\eta^{\p}\|_{\infty}^2\|w\|_{\infty}^2-C_4 e^{-C\wt}\big]k^2=
0.
\end{eqnarray}
%

Set the constant

%
\begin{eqnarray}\label{EQconstant5}
C_5=\frac{1}{2}\left(1+[\frac{1}{2 M^2}-1]_{+}\right)\sup_{t\in[0,T)}\|u(t)\|^2_{\W^{1,\infty}},
\end{eqnarray}
%

where $[\cdot]_{+}$ denote the positive part of a real number; therefore by
(\ref{EQnonlinear8}), (\ref{EQab}), (\ref{EQxy}), the Lipschtiz regularity of $u$, and using the Young inequality with exponent $2$ we can write,

%
\begin{eqnarray}\label{EQnonlinear9}
\lefteqn{\Delta^2 G(v_1 e^{C \wt},v_2 e^{C \wt},2p+m,2q+m)=}\nonumber\\
&=& \Delta^2 (G- C_4 e^{-C \wt}\frac{u^{2}}{2})(v_1 e^{C \wt},v_1 e^{C \wt},2p+m,2q+m)+\nonumber\\
&&+\frac{1}{4} C_4 e^{C \wt}|v^2_2-v^2_1|^2\leq \nonumber\\
&\leq&\frac{1}{4} C_4 e^{C \wt}|v_1-v_2|^2=\nonumber\\ 
&=&\frac{1}{4} C_4 e^{-C\wt}|u_1-u_2|^2\leq\nonumber\\
&\leq& \frac{1}{2} C_4 \left(|u_1-u_3|^2+|u_2-u_3|^2\right)\leq\nonumber\\
&\leq& \frac{1}{2} C_4 \sup_{t\in[0,T)}\|u(t)\|^2_{\W^{1,\infty}}
[(\frac{\delta |p|}{2 M})^{\frac{2}{3}}+(\frac{\delta |q|}{2 M})^{\frac{2}{3}}]
\leq\nonumber\\
&\leq& \frac{1}{2}C_4 \sup_{t\in[0,T)}\|u(t)\|^2_{\W^{1,\infty}}\left(\frac{1}{2}[(\frac{\delta}{2 M})^{\frac{1}{3}}(|p|^{\frac{4}{3}}+|q|^{\frac{4}{3}})]+\frac{\delta}{2 M}\right)\leq\nonumber\\ 
&\leq& \frac{1}{2}C_4 \sup_{t\in[0,T)}\|u(t)\|^2_{\W^{1,\infty}}\left((2v_3-v_1-v_2)+\delta M [\frac{1}{2 M^2}-1]_{+}\right)\leq\nonumber\\ 
&\leq& C_4 C_5 (2v_3-v_1-v_2).
\end{eqnarray}


Introducing estimates (\ref{EQnonlinear2}), (\ref{EQnonlinear3}), (\ref{EQnonlinear4}), 
(\ref{EQnonlinear6}), (\ref{EQnonlinear9}), in (\ref{EQnonlinear1}), we finally obtain

%
\begin{eqnarray}\label{EQnonlinearlast}
\lefteqn{\big[\l(v_1 e^{C \wt})e^{C \wt }|\s^{\top}\Phi_1|^2
-\eta (v_1 e^{C \wt})\la \s^{\top} \Phi_1, w_1 \ra+\l(v_2 e^{C \wt})e^{C \wt }|\s^{\top}\Phi_2|^2}\nonumber\\
&&-\eta (v_2 e^{C \wt})\la \s^{\top} \Phi_2, w_2 \ra
-2\l(v_3 e^{C \wt})e^{C \wt}|\frac{1}{2}\s^{\top}\Phi_3|^2
-\eta(v_3 e^{C \wt})\la \s^{\top}\Phi_3,w_3\ra\big]
\leq \nonumber\\
&&\leq(C_1+C_2+C_3+C_4C_5)(2v_3-v_1-v_2)+O(\sqrt{\var}).
\end{eqnarray}


Consider now the last term in (\ref{EQadding}), by the regularity assumptions on the function $f$, the same
argument of (\ref{EQnonlinear9}) used for estimating $|u_1-u_2|$ and again (\ref{EQxy}), we have,

%
\begin{eqnarray}\label{EQsource1}
f_1+f_2-2 f_3 &\leq& 
\sup_{t\in [0,t)}\|f(t)\|_{\W^{2,\infty}}\big[\big[(\frac{\delta}{2 M})^{\frac{4}{3}}(|p|^{\frac{4}{3}}+|q|^{\frac{4}{3}})+
(\frac{\delta}{2 M})^2 |m|^2]^{\frac{1}{2}}+\nonumber\\
&&+|u_1-u_3|^2+|u_2-u_3|^2+|u_1+u_2-2 u_3|\big]\leq \nonumber\\
&\leq& \sup_{t\in [0,t)}\|f(t)\|_{\W^{2,\infty}} \big\{2v_3-v_1-v_2+\delta M[\frac{1}{4 M^2}-1]_{+}\nonumber\\
&&+2 \sup_{t\in[0,T)}\|u(t)\|^2_{\W^{1,\infty}}C_5 (2 v_3-v_1-v_2)+e^{C \wt}(2v_3-v_1-v_2)\big\} 
\leq \nonumber\\
&\leq& C_6(2v_3-v_1-v_3),
\end{eqnarray}
%

where we consider

%
\begin{eqnarray}\label{EQconstant6}
C_6=\sup_{t\in [0,t)}\|f(t)\|_{\W^{2,\infty}}\left(1+[\frac{1}{4 M^2}-1]_{+}+2 \sup_{t\in[0,T)}\|u(t)\|^2_{\W^{1,\infty}}C_5
+e^{C\wt}\right).
\end{eqnarray}
%

Finally using (\ref{EQlinear22}), (\ref{EQlinear1}), (\ref{EQnonlinearlast}), and (\ref{EQsource1}) in (\ref{EQadding}), we have

%
\begin{eqnarray}\label{EQadding1}
\frac{\gamma}{T^2}&\leq& [6 \sup_{t\in (0,T)}\|\mu(t)\|_{\W^{2,\infty}}+C_1+C_2+C_3+C_4C_5\nonumber\\
&&+C_6 e^{-C\wt}-C](2v_3-v_1-v_3)+O(\sqrt{\var}).
\end{eqnarray}


From the definition (\ref{EQconstant6}), we see that $C_6 e^{-C\wt}$ is bounded as a function of $C>0$, so
if we choose $C$ sufficently great, we obtain a contradiction letting $\var\rightarrow 0$. This prove the result.
\cvd
\end{Dimo}


Now we use these Propositions to eliminate the conditions on $\l$ and to obtain the assertion
of the Theorems \ref{SEMICONVEX} and \ref{SEMICONCAVE}.\\

\begin{Dimo} {\bf Theorem \ref{SEMICONVEX}.}
We go to build a change of variable such that the new differential equation will have
the required structural properties of Proposition \ref{regularity}.\\
Let $c,\L$ be respectively the infimum of the closed interval $I$ where $u$ take its values, and the primitive 
of $\l$ with $\L(c)=0$. Then consider the solution $Q=Q(\tau)$ of the following ordinary Cauchy problem
\begin{eqnarray}\label{CHG}
\left\{\begin{array}{c}
  \frac{d Q}{d \tau}= \exp( 4 \sqrt{\tau+1}+ 2\L(Q)), \\
Q(0)=c.
\end{array}\right.
\end{eqnarray}
The problem admits an increasing local solution. Moreover
by the continuity of $\l$, $Q$ is also $\C^2$. We prove that $Q$ maps the interval $I$. Consider the following cases
\begin{description}
\item[{\bf Case 1, ($b=\infty$)}.] Let $[0,\tau^{\star})$ be the maximal interval of existence for $Q$, and denote with
$\o{Q}$ the limit for $\tau\rightarrow \infty$.
If $\tau^{\star}=\infty$, then $\o{Q}=\infty$. Actually by the equation (\ref{CHG}), we have, 
\begin{eqnarray}\label{CHG1}
Q^{\p}(\tau)\geq \exp\big(\inf_{u\in[c,\o{Q}]}\L(u)+4\big),\;\;\;\forall\;\tau>0.
\end{eqnarray} 

So integrating (\ref{CHG1}) from $0$ to $\tau>0$ and letting $\tau$ to infinity we obtain the assertion. 
If $\tau^{\star}<\infty$, then for definition of maximal interval $Q$ blow-ups at $\tau^{\star}$.
Otherwise since that $\l$ is defined in $(a,\infty)$, the solution $Q$ could be extended.\\ 
\item[{\bf Case 2, ($b<\infty$)}.] Consider again the maximal interval of existence. Then
if $\tau^{\star}=\infty$ and $\o{Q}$ is strictly less than $b$, then, again using (\ref{CHG1}), we obtain a contradition.
If $\tau^{\star}<\infty$, then and $\o{Q}<b$, then the solution $Q$ can be extended bacuase
$\l$ is continuous in $(\o{Q},b)$.
\end{description}

In each case the function reachs $b$ in the limit sense; in particular $Q$ maps $I$. Moreover
$Q$ can be defined in an open interval $V\subset (-\var_0,\infty)$, $\var_0>0$, and
$I\subset Q(V)\subset (a,b)$. By the increasing property of $Q$, it admits a $\C^2(Q(V);V)$ inverse, which we donote
$P$. We use $Q$ as a transformation for a global change of the variable $u$ . The function $\tau=P\circ u$, is a bounded, $t$-uniformly Lipschitz continuous viscosity solution of
\begin{eqnarray}\label{CHG2}
\partial_t \tau-\frac{1}{2}tr(\s\s^{\top} \nabla^2 \tau)+\la \mu, \nabla \tau \ra-\frac{|\s^{\top} \nabla \tau|^2}{\sqrt{\tau+1}}\nonumber\\
+\eta(Q(\tau))\la \s^{\top} \nabla\tau, w \ra +f(x,t, Q(\tau))=0,\;\;\;(x,t)\in\R^N\times(0,T),
\end{eqnarray}
where the initial datum is $\tau_0=P\circ u_0$, and which takes values in the closed interval $P(I)$. It is easy to verify the structural hypothesis $i)$, $ii)$
of Proposition \ref{regularity}, where $\l$ and $\eta$ are subsituted, respectively, by the functions $-(1+\tau)^{-\frac{1}{2}}$ and $\eta\circ Q$, over the interval $V$, with regularity
properties over $P(I)$. So applying Proposition \ref{regularity}, we deduce that there exist positive constants $C,K_0>0$, such that, 
\begin{eqnarray}\label{CHG3}
\tau(x+h,t)+\tau(x-h,t)-2\tau(x,t)\geq -e^{C t} K_0 |h|^2,\;\;\forall\;x,h\in\R^N,\;t\in[0,t),
\end{eqnarray}
where the constant $K_0$ depends on $P$ and, on $L_0$ and $Lip(u_0)$. Therefore for every $x,\;h\in \R^N$, $t\in[0,t)$, and some $s^+,\;s^-\in [0,1]$, denoting
\begin{eqnarray}\label{CHG5}
\tau^+&=&\tau(x,t)+s^+(\tau(x+h,t)-\tau(x,t)),\nonumber\\
\tau^-&=&\tau(x,t)+s^-(\tau(x-h,t)-\tau(x,t)),\nonumber\\
\tau&=&\tau(x,t),
\end{eqnarray}
we can write,
\begin{eqnarray}\label{CHG4}
&&u(x+h,t)+u(x-h,t)-2u(x,t)=Q^{\p}(\tau^+)(\tau(x+h,t)-\tau(x,t))\nonumber\\
&&+Q^{\p}(\tau^-)(\tau(x-h,t)-\tau(x,t))\nonumber\\
&=& Q^{\p}(\tau)(\tau(x+h,t)+\tau(x-h,t)-2\tau(x,t))+s^+ Q^{\p\p}(\tau^{++})(\tau(x+h,t)\nonumber\\
&&-\tau^2(x,t)+ s^- Q^{\p\p}(\tau^{--})(\tau(x-h,t)-\tau^2(x,t)\nonumber\\
&&\geq -[K_0 C_{00}+C_0]|h|^2.
\end{eqnarray}
for some $\tau^{++}\in [\min (\tau,\tau^+),\max (\tau,\tau^+)]$, $\tau^{--}\in [\min (\tau,\tau^-),\max (\tau,\tau^-)]$, where 
$C_{00}$ is a positive constant depending on $Q^{\p}$, while $C_0$ depends on the Lipschitz constant of the solution $u$ and on $Q^{\p\p}$.
This prove the assertion of the Theorem.
\cvd
\end{Dimo}

The equivalent result for the semiconcavity property can be then obtained with same arguments.
So we limit us to give some outlines in the following proof.\\

\begin{Dimo} {\bf Theorem \ref{SEMICONCAVE}.}
By the same notations used for proving Theorem \ref{SEMICONVEX}, we observe that, choosing
the increasing transformation $u=Q(\tau)$, where $Q$, is implicitly defined as the solution of the 
oridnary Cauchy problem,
\begin{eqnarray}\label{CHG5}
\left\{\begin{array}{c}
 \frac{d Q}{d \tau}=\exp( -\frac{2}{l+1}(\tau+1)^{l+1}+ 2\L(Q)), \\
Q(0)=c,
\end{array}\right.
\end{eqnarray}
where $l$, is choosen as bigger than $3$, then we obtain, as in the semiconvexity case a new equation
for $\tau$, which satisfies the structural hipothesis required for applying Proposition \ref{regn}.
\cvd
\end{Dimo}

As an immediate consequence of Theorems \ref{SEMICONVEX} and \ref{SEMICONCAVE},
and by the definition of the space $\W^{2,\infty}$, it follows the proof of
the Theorem \ref{regf1}. This conclusive fact allows us to have
a second-order regularity result for the solution of problem (\ref{EQ}), with regular initial data. 

\section{Time regularity}
In this section we use the spatial regularity of the solution $u$, obtained through the Theorem \ref{regf1}, to prove the time regularity of it.\\ 
The result which we show here, see Theorem \ref{TIME}, is not stated for second order Hamilton-Jacoby equations which have the structure of (\ref{EQ}), since the lack of regularity in the spatial 
variable. So we present it as a possible interesting extension of previous works, in the framework of the viscosity theory.\\ 
In order to simplify the notations, in the following, for every function $g$ defined in 
$\R^N\times[0,T)\times(a,b)$, and $T>h>0$, we set $g_h(x,t,u)=g(x,t+h,u)$, 
where $0\leq t<T-h$, and $x\in\R^N$, $u\in(a,b)$. Moreover, if 
$g$ is a Lipschitz continuous function over $[0,T)$, uniformly with respect to the other variables, then we shall 
denote $L(g)$ the constant defined as:
\begin{eqnarray*}
L(g)=\sup_{\begin{array}{c}
t,s\in[0,T)\\ 
t\neq s\\
(x,u)\in\R^N\times(a,b)
\end{array}
}\frac{|g(x,t,u)-g(x,s,u)|}{|t-s|}
\end{eqnarray*} 
\begin{teo}\label{TIME}
Let $u(t)\in\W^{2,\infty}(\R^N)$, for every $t\in[0,T)$. Assume that $\s,\;\mu,\;w,\;f$ are Lipschitz continuous functions of the time, uniformly with
respect to the other variables, and $\mu,\;w,\;f$ are bounded, then $u(x,\cdot)$, is a Lipschitz continuous function uniformly in $\R^N$.
\end{teo}

This result is based on the next Lemma, which shows how the solution $u$ tends to the initial datum, when we send the time at zero.

\begin{lem}\label{U0DIP}
If $u_0\in\W^{2,\infty}$, and $\mu,\;w,\;f$ are bounded, then there exists
a constant $C_0>0$, such that
\begin{eqnarray}\label{U0DIP1}
|u(x,t)-u_0(x)|\leq C_0 t,
\end{eqnarray}
holds for every $(x,t)\in \R^N\times [0,T)$.
\end{lem}

\begin{Dim}.
We limit us to prove the inequality,
\begin{eqnarray}\label{U0DIP2}
u-u_0\geq -C_0 t,
\end{eqnarray}
choosing,
\begin{eqnarray}\label{U0DIP7}
C_0\!\!&=&\!\!\sup\{H(x,t,u,p,X)\;:\;(x,t,u)\in\R^N\times [0,T)\times I,\;|p|\leq\|\nabla u_0\|_{\infty},\nonumber\\
&& \|X\|\leq\|\nabla^2 u_0\|_{\infty}\}.
\end{eqnarray}
The other follows with same arguments.
Consider the function,
\begin{eqnarray}\label{U0DIP3}
\Phi_{\gamma,\var}(x,t)=u(x,t)-u_0^{\var}(x)+C_0 t+\var|x|^2+\frac{\gamma}{T-t},
\end{eqnarray}
For every $\gamma,\var>0$, $(x,t)\in\R^N\times[0,T)$. Where $u_0^{\var}$ is the convolution between the
standard mollifier and $u_0$. The inequality (\ref{U0DIP2}) is satisfied if for every $\gamma>0$,
there is $\var_0>0$, such that, for every $0<\var<\var_0$, 
\begin{eqnarray}\label{U0DIP4}
\inf_{\R^N\times[0,T)}\Phi_{\gamma,\var}\geq 0.
\end{eqnarray}
Actually if (\ref{U0DIP4}) holds, setting a point $(y,s)\in\R^N\times(0,T)$, we have
\begin{eqnarray*}
u(y,s)-u_0^{\var}(y)+C_0 s+\var|y|^2+\frac{\gamma}{T-s}\geq 0.
\end{eqnarray*}
So we send $\gamma,\;\var$ to zero and obtain (\ref{U0DIP2}).\\
Suppose that (\ref{U0DIP4}) would be false. Hence by the same observations used for proving the Proposition
\ref{regularity}, there is $\gamma_0>0$, and a sequence $\var_j\rightarrow 0$, as $j\rightarrow \infty$,
such that $\Phi_{\gamma,\var}$ has a global minimum point, $(\wx,\wt)$, and $\Phi_{\gamma,\var}(\wx,\wt)<0$,
for $\gamma=\gamma_0$, $\var=\var_j$. If $\wt=0$, then
\begin{eqnarray*}
0>\Phi_{\gamma,\var}(\wx,0)\geq -\var \|u_0\|_{\W^{1,\infty}}+\frac{\gamma}{T}.
\end{eqnarray*}
Hence for large $j$, $\wt>0$, and the minimum point is an interior point; moreover, $u_0^{\var}-C_0 t-\var|x|^2-\frac{\gamma}{T-t}$ is a test function for $\P^{2,-}u$ at $(\wx,\wt)$.
By (\ref{U0DIP3}), we deduce that
\begin{eqnarray}\label{U0DIP5}
\var|\wx|=\sqrt{\var^2|\wx|^2}= O( \var^{\frac{1}{2}}).
\end{eqnarray}
Moreover by the regularity assumption on $u_0$, we have $\|\nabla u_0^{\var}\|_{\infty}\leq \|\nabla u_0\|_{\infty}$ 
and $\|\nabla^2 u_0^{\var}\|_{\infty}\leq \|\nabla^2 u_0\|_{\infty}$, where $\nabla u_0,\;\nabla^2 u_0$, denote
the weakly derivatives of $u_0$. Introducing the derivatives of the test function in the equation (\ref{EQ}) and using the assumption
of boundness of the coefficients, we can write
\begin{eqnarray}\label{U0DIP6}
C_0+\frac{\gamma}{T^2}\leq H(\wx,\wt, u(\wx,\wt), \nabla u_0^{\var}(\wx,\wt),\nabla^2  u_0^{\var}(\wx,\wt))+O(\var^{\frac{1}{2}}).
\end{eqnarray}
By (\ref{U0DIP7}), letting $j\rightarrow \infty$ in (\ref{U0DIP6}), we obtain the contradiction.
\end{Dim}
\cvd

\begin{Dimo} {\bf Theorem \ref{TIME}.}
Consider $u_{h}$, for $T>h>0$. $u_h(t)\in\W^{2,\infty}(\R^N)$ uniformly in time $t<T-h$, and is a viscosity solution of the problem

\begin{eqnarray}\label{TIME1}
&&\partial_t u_h-\frac{1}{2}tr(\s\s^{\top}(t+h)\nabla^2 u_h)+\la\mu_h ,\nabla u_h \ra+\l(u_h)|\s^{\top}(t+h)\nabla u_h|^2\nonumber\\
&&\;\;\;+\eta(u_h)\la \s^{\top}(t+h)\nabla u_h,w_h\ra +f_h(x,t,u_h)=0,
\end{eqnarray} 

in $\R^N\times (0,T-h)$. For every $(x,t)\in\R^N\times[0,T)$, define

\begin{eqnarray}\label{TIME2}
u^h(x,t)&=&u(x,t)+(C_0+\alpha(t))h,\nonumber\\
\alpha(t)&=& \frac{e^{B_1 t}-1}{B_1}( C_0 B_1 +B_2).
\end{eqnarray}

Where $C_0$ has been defined in Lemma \ref{U0DIP}, while
\begin{eqnarray}\label{TIME3}
B_1 &=& \|\l\|_{\infty}\|\s^{\top}\|^2_{\infty}\sup_{t\in [0,T)}\|u(t)\|_{\W^{1,\infty}}^2\nonumber\\
&&+\|\eta\|_{\infty}\|\s^{\top}\|_{\infty}\|w\|_{\infty}\sup_{t\in [0,T)}\|u(t)\|_{\W^{1,\infty}}+
L(f).\nonumber\\
\\
B_2 &=& L(\s\s^{\top})\sup_{t\in [0,T)}\big(\frac{1}{2}N^2 \|u(t)\|_{\W^{2,\infty}}+\|\l\|_{\infty}
\sup_{t\in [0,T)}\|u(t)\|^2_{\W^{1,\infty}}\big)\nonumber\\
&&+\sup_{t\in [0,T)}\|u(t)\|_{\W^{1,\infty}}\big(L(\s^{\top})\|w\|_{\infty}+L(w)\|\s^{\top}\|_{\infty}+L(\mu)\big).\nonumber
\end{eqnarray}
Then $u^h$ is a viscosity solution of the following problem
\begin{eqnarray}\label{TIME4}
&&\partial_t u^h-\alpha^{\p}h-\frac{1}{2}tr(\s\s^{\top}\nabla^2 u^h)+\la \mu, \nabla u^h \ra+
\l(u^h-C_0 h-\alpha h)|\s^{\top}\nabla u_h|^2\nonumber\\
&&\;\;\;+\eta(u^h-C_0 h-\alpha h)\la \s^{\top}\nabla u^h,w\ra +f(x,t,u^h-C_0 h-\alpha h)=0.
\end{eqnarray}
Recall the Hamiltonian $H$ given by (\ref{H}), then





we observe that by the regularity assumption on $u(t)$ for every time $t$, we have, for every element $(s,p,X)\in\P^{2,\pm}(x_0,t_0)$, for $(x_0,t_0)\in\R^N\times(0,T)$,
$|p|\leq \|u(t_0)\|_{W^{1,\infty}}$, $\|X\|\leq C \|u(t_0)\|_{W^{2,\infty}}$, for some constant $C>0$. \\












Therefore for every $(x,t)\in\R^N\times(0,T-h)$, if $(s,p,X)\in\P^{2,+}u_h(x,t)$, then
$(s,p,X)\in\P^{2,+}u(x,t+h)$, and  $|p|\leq \sup_{t\in[0,T)}\|u(t)\|_{\W^{1,\infty}}$ while
$|X|\leq C \sup_{t\in[0,T)}\|u(t)\|_{\W^{2,\infty}}$. The same kind of observation holds
for $u^h$.\\
Moreover using the last observations, and introducing the constants (\ref{TIME3}), we can estimate (\ref{TIME1}) and (\ref{TIME4}), obtaining
\begin{eqnarray}\label{TIME6}
&&\partial_t u_h+H(x,t,u_h,\nabla u_h,\nabla^2 u_h)-B_2 h\leq 0,\nonumber\\
&&\partial_t u^h+H(x,t,u^h,\nabla u^h,\nabla^2 u^h)+B_1(\alpha+ C_0)h-\alpha^{\p}h\geq 0,
\end{eqnarray}
on the domain $\R^N\times (0,T-h)$, in a viscosity sense. By the definition (\ref{TIME2}), we have
\begin{eqnarray}\label{TIME7}
0 &\leq& \partial_t u^h+H(x,t,u^h,\nabla u^h,\nabla^2 u^h)+B_1(\alpha + C_0)h-\alpha^{\p}h\nonumber\\
 &= & \partial_t u^h+H(x,t,u^h,\nabla u^h,\nabla^2 u^h)-B_2 h.
\end{eqnarray}
Hence $u_h$ and $u^h$ are subsolution and supersolution of the same problem, repsectively, and by Lemma
\ref{U0DIP}, $u_h(x,0)=u(x,h)\leq u_0+C_0 h=u^h(x,0)$, for every $x\in\R^N$. Using the comparison principle, which
holds applying the comparison principle which we have stated in \cite{38}, we deduce that $u(x,t+h)-u(x,t)\leq C_0 h+\alpha(t)h$, for every $(x,t)\in\R^N\times [0,T)$.
Using the same argument, with the same function $\alpha$, and using again Lemma \ref{U0DIP}, we have
$u(x,t+h)\geq u(x,t)-C_0 -\alpha(t)h$, for each point $(x,t)\in \R^N\times [0,T-h)$. These two inequalities imply
that, for every $x\in\R^N$, $t_1,\;t_2\in[0,T)$, 
\begin{eqnarray}\label{TIME8}
|u(x,s_2)-u(x,s_1)|\leq C_T|s_2-s_1|.
\end{eqnarray}
Where the constant $C_T$ depends on the final time $T$, and $B_1,\;B_2\;C_0$.
\end{Dimo}
\cvd

By Theorem \ref{regf1} and Theorem \ref{TIME}, then it easily follows Theorem 
\ref{SPACETIME}.
 
\section{Regularity and the Ito's formula}\label{Conclusions}
In this section, we derive some consequences from the results of the previous sections. Let us recall first from \cite{KOCAN} a useful result about 
Sobolev spaces $\W^{2,1,p}_{loc}(\R^N\times(0,T))$, i.e.: the space of the functions $u$, such that $u,\;\nabla u,\;\nabla^2 u,\;\partial_t u$ $\in L^{p}_{loc} (\R^N\times (0,T))$.\\ 
We start by an observation about the Lebesgue points.

\begin{oss}\rm\label{2.2}
If $g\in L^p_{loc}(\R^{N+1})$, for every $p\in[1,\infty)$, then for a.e. $(\wx,\wt)$, 
\begin{eqnarray*}
\lim_{r\rightarrow 0^+} \frac{1}{|\Omega_r(\wx,\wt)|}\int_{\Omega_r(\wx,\wt)} |f-f(\wx,\wt)|^p dx dt=0.
\end{eqnarray*}
where $\Omega_r(x,t)=B_r(x)\times (t-r^2,t)$, and $B_r(x)$ denotes the open ball in $\R^N$ with ray $r$, and centre at $x$.
This follows from a version of Vitali's covering Theorem by replacing the balls with $\Omega_r$, see e.g.
Remark I.3.1 in \cite{GUZMAN}.
\end{oss}

\begin{prop}\label{KOCAN1}
Let $p>\frac{N+2}{2}$, $N\geq 2$, and $u\in\W^{2,1,p}_{loc}(\R^N\times (0,T))$. Let $(\wx,\wt)$ be a Lebesgue point
(in the sense of (\ref{2.2})) of $u$, $\nabla u$ (in $L^1$), and $\nabla^2 u$, $\partial_t u$ (in $L^p$).
Then as $(x,t)\rightarrow (\wx,\wt)$, we have
\begin{eqnarray}\label{KOCAN2}
&& |u(x,t)-\big( u(\wx,\wt)+\partial_t u(\wx,\wt)(t-\wt)+\la\nabla u(\wx,\wt),x-\wx\ra\nonumber\\
&& \;\;\;+\frac{1}{2}
\la \nabla^2 u(\wx,\wt)(x-\wx),(x-\wx)\ra\big)|\leq o(|x-\wx|^2+|t-\wt|).
\end{eqnarray}
\end{prop}
For the proof we refer the reader to \cite{KOCAN}.
If the assumptions of Theorem \ref{SPACETIME} are satisfied, then we can apply Proposition \ref{KOCAN1} to $u$, obtaining that for a.e. $(x,t)\in\R^N\times (0,T)$, we have 
\begin{eqnarray}\label{KOCAN2}
\partial_t u(x,t)+H(x,t,u(x,t),\nabla u(x,t),\nabla^2 u(x,t))=0,
\end{eqnarray} 
where $\partial_t u, \nabla u,\nabla^2 u$, represent the wealkly derivatives of $u$.
This regularity of the solution $u$ has a useful application when we use the Ito's rule along a
stochastic process.\\
As we have noted in the Introduction, the representation of financial derivatives can be expressed
as a conditional expectation over the probability space of the underlying factors that determine the instrument's price.\\
In the equality (\ref{GABAIX}), we introduce a relation between this conditional expected value and the solution of (\ref{DM1}). This 
is a consequence of the regularity of $u=U+h+\xi$, which follows by applying the results of the previous sections to (\ref{DM1}). In the model proposed by X. Gabaix in \cite{19}, 
the measure $Q$ depends in a nonlinear way by the price of the Mortgage-Backed Security and its volatility; this feature produces the quadratic 
nonlinear term in the equation (\ref{DM1}). To rigorously derive the equation (\ref{GABAIX}), we need to apply
the Ito's formula to $u(X_t,t)$, where $X_t$ is the underlying stochastic process. When $u$ is less regular 
than $\C^{2,1}$, a first result was established by Krylov \cite{KR}. In a more recent work, 
Haussmann \cite{ITO} describes a result in this direction, which
shows that the Ito's formula holds for $u\in\W^{2,1,\infty}(\R^N\times(0,T))$, provided that it is interpreted appropriately,
using the generalized Hessian.\\ 
We use the regularity of the solution $U$, which gives the property (\ref{KOCAN2}), combining it with the result of Haussmann to obtain the equality (\ref{GABAIX}).

In the next we shall assume $\mu\in\W^{2,1,\infty}(\R^N\times[0,T]$ and for sake of simplicity we limit us to treat in detail the case of a constant matrix $\s$ and consider the case $N\geq d$.\\
As we have already pointed out in the Introduction, equation (\ref{DM1}) is 
equivalent to the general problem (\ref{EQ}). Actually, after the change $u=U+h+\xi$, $u$ solves a problem which has the same structure of (\ref{EQ}), where in particular $w=\s^{\top}\cdot\nabla h$.\\
Hence, following the comparison principle and the Lipschitz regularity for viscosity solutions of (\ref{DM1}) proved in \cite{38}, and assuming $h$ to be a smooth function, as in \cite{38},  by Theorem \ref{SPACETIME} is straightforward to assume the existence of a unique viscosity 
solution $U\in\W^{2,1,\infty}(\R^N\times[0,T))$ of the equation (\ref{DM1}) with $U(x,0)\equiv 0$, such that $U+h+\xi>0$.

\begin{teo}\label{HAUSSMANN}
Let  $(\Omega,\F,P)$ be a probability space and let $(W_t,\F_t)_{t\in [0,T]}$ be a $d$-dimensional continuous standard Brownian motion
over $\Omega$ and let $X_t$ be the Ito process defined as the solution of the stochastic differential equation (\ref{GABAIX.1}), for some intial datum $X_0\in\R^N$. Define 
\begin{eqnarray}\label{ITO1}
\gamma_s &=& \rho\frac{\s^{\top}(T-s)\nabla U(X_s,T-s)}{U(X_s,T-s)+h(X_s,T-s)+\xi(T-s)},\;\;\mbox{a.s.}\nonumber\\
\widehat{W}_s&=&W_s+\int_0^s \gamma(\kappa)d\kappa,\nonumber\\
\frac{dQ}{dP}&=&e^{-\int_0^T \gamma_s ds-\frac{1}{2}\int_0^T |\gamma_s|^2 ds},
\end{eqnarray} 
with $0<\rho<1$ a parameter. If $N=d$ we do not need any assumption.\\
Otherwise for the case $N>d$, assume that the {\em projection-process} $\pi_t$, defined as the projection of $X_t$ over the Kernel of $\s^{\top}$, satisfies
\begin{eqnarray}\label{N>d}
A\mapsto P\left(\pi_s\in A| X_t \right),\;A\subset R^m\;\mbox{{\em is ${\cal L}^m$-absolutely continuous,}} 
\end{eqnarray}
for every $s>t$, where ${\cal L}^m$ denotes the Lebesgue masure over $\R^m$ and $m=N-rank(\s^{\top})$. Then 
\begin{eqnarray}\label{GABAIX1}
U(X_t,T-t)=\E_t^{Q}\Big[\int_t^T \big(\tau-r(T-s)\big) e^{-\int_t^s r(T-\kappa)d\kappa} h(X_s,T-s) ds\Big],
\end{eqnarray}
holds $a.s.$ for every $0<t\leq T$.
\end{teo}

We remark that a classical result, the Girsanov's Theorem, state that $\widehat{W}_s$ is a Brownian motion with the same filtration of $W_s$, see for istance \cite{10} or \cite{9}.
To prove the relation (\ref{GABAIX1}) for the solution $U$ to the model (\ref{DM1}), we shall use the results of \cite{ITO}. Actually the assertion of Haussmann interprets the Ito's rule through some processes which substitute the usual derivatives of the function evaluated at the process $(X_t,t)$, which are equal if the underlying process belongs to some set of full Lebesgue measure.
In order to state our formula we need to neglect the term which correponds
to an integration over the paths of $X$ which fall in a set of null measure. 
Therefore the assumption (\ref{N>d}) can be explained by that purpose. Actually the degeneration of the volatility matrix plays a crucial role: in general along the components of the vector $X$ which correspond to the directions of the kernel of $\s^{\top}$, we have {\em almost-deterministic} trajectories which surely form a set of zero measure; this feature, in general, could not allow to justify the previous assertion (\ref{GABAIX1}) about $U$, using the Haussmann's Theorem, unless the drift coefficient has some compatibility condition with $\s$; in fact this feature is expressed through the condition on the projection-process $\pi_t$, see the next remark \ref{osservazione}.\\
As we have pointed out in the Introduction, there are some works such as \cite{BELL}, where the existence of densities for the solutions of stochastic differential equations under the H$\ddot o$rmander's condition is proved. In particular Theorem 3.1 of chapter 3 by \cite{BELL} gives a sufficient condition, which involves the use of the Lie bracket between the set of vector fields represented by the drift coefficient and the columns of the volatility matrix, to prove the existence of an absolutely continuous distribution with respect to the Lebesgue measure. Although that condition can be used in practice to investigate differential equations with a complicated form, it does not allow to consider a time dependence of the coefficients, and it requires also more regularity for computing the Lie brackets of the vector fields. Moreover, while for financial purposes the coefficients $\mu$ and $\s$ can be taken as constants, the probabilistic H$\ddot o$rmander's condition used by D.R. Bell in \cite{BELL} does not hold in this case. Therefore we propose an alternative assumption, contained in (\ref{N>d}), which covers the cases of interest for our applications.  Actually let us notice that case $N=d$, the assumption (\ref{N>d}) can be removed via an approximation of the underlying process $X$. The proof of the following technical Lemma motivates of condition (\ref{N>d}).  

\begin{lem}\label{HAUSSMANN0}
Let $Q$ be a probability measure equivalent to $P$ over $\Omega$ and let $B$ be a set of zero Lebesgue measure over $\R^N\times[t,T]$, for $0\leq t<T$. Set 
\begin{eqnarray}\label{HAUSSMANN0.1}
Z_{t,T}=\int_t^T p_s 1_B(X_s,s)ds+q_s 1_B(X_s,s) dW_s,
\end{eqnarray}
where $p_s,q_s$ are $(\F_{s})_{t\leq s\leq T}$-adapted bounded processes. Then, under the assumption (\ref{N>d}), it holds
\begin{eqnarray*}
\E_t^{Q}\big[|Z_{t,T}|\big]=0.
\end{eqnarray*}
\end{lem}
 
\begin{oss}\label{REMARK}
Let us recall the following important fact. For every integrable random variable $\F_T$-measurable there holds
\begin{eqnarray*}
\E_t^Q\left[X\right]=\frac{\E_t\left[X \frac{dQ}{dP}\right]}{\E_t\left[\frac{dQ}{dP}\right]},
\end{eqnarray*}
where $\E_t[\cdot]$ denote the conditional mean with respect to the measure $P$.
To see  this relation between the conditional means taken through two equivalent probability measure, we refer the reader to \cite{13}. 
\end{oss}

\begin{Dimo} {\bf Lemma \ref{HAUSSMANN0}.}
By remark \ref{REMARK}, the  boundness of the processes $p_s$, $q_s$ and 
\[\E_t\left[\int_t^T q_s 1_B(X_s,s)dW_s\right]\leq \E\left[\left|\int_t^T q_s 1_B(X_s,s)dW_s\right|^2\right]^{\frac{1}{2}}\]
\[=E\left[\int_t^T |q_s|^2 1_B(X_s,s)ds\right]^{\frac{1}{2}},\]
where in the last passage we have used the classical Ito isometry, we are reduced to prove
\begin{eqnarray}\label{HAUS0}
\E_t\left[\int_t^T 1_{B}(X_s,s)ds\right]=0.
\end{eqnarray}
{\bf Case $N>d$.} Let $\{b_1,\ldots,b_m\}$ be an orthonormal basis of the kernel of $\s^{\top}$, then consider an $N\times N$ invertible matrix $M$, such that
\[M^{\top}e_i=b_i,\qquad i=1,\ldots,m,\] 
and define the process $Y_s=M X_s$. Therefore 
\[Y_s=(\pi_s^M,G_s):=\left((M^{\top})^{-1}_{|_{Ker(\s^{\top})}}\pi_s,G_s\right),\]
and
\begin{eqnarray}\label{HAUS1}
d \pi^M_s &=& \mu_{\pi}(Y_s,s)ds,\qquad \pi^M_s\in\R^m,\nonumber\\
\\
dG_s &=& \mu_G(Y_s,s)ds+\s_G\cdot dW_s,\qquad G_s\in\R^{N-m}.\nonumber
\end{eqnarray}
Moreover
\begin{eqnarray}\label{HAUS2}
\mu_{\pi}(y,s)&=&(\la \mu(M^{-1}y,T-s),b_1\ra,\ldots,\la\mu(M^{-1}y,s),b_m\ra),
\end{eqnarray}
and $\s_G \in \M_{N-m,d}(\R)$ is a matrix of rank $N-m$. Let 
\begin{eqnarray*}
\t{M}=\left(\begin{array}{ccc}
  &    & 0 \\
  &  M & \vdots \\
  &    & 0\\
0 & \ldots & 1
\end{array}\right),\qquad \t{B}=\t{M}B.
\end{eqnarray*}
Noting that $\t{B}$ has zero Lebesgue measure, (\ref{HAUS0}) is equivalent to
\begin{eqnarray}\label{HAUS3.7}
\E_t\left[\int_t^T 1_{\t{B}}(Y_s,s)ds\right]=0.
\end{eqnarray}
Fix $\var>0$ and $t^+>t/(1-\var)$, then consider the approximating process $Y_s^{\var}=(\pi_s^{M,\var},G_s^{\var})$, where 
\begin{eqnarray}\label{HAUS4}
\pi_s^{M,\var}&=&\pi^M_t+\int_t^{s(1-\var)}\mu_{\pi}(Y_{\l},\l)d\l\nonumber\\
\\
G_s^{\var}&=&G_t+\int_t^{s(1-\var)}\mu_G(Y_{\l},\l)d\l+\s_G\cdot(W_s-W_t).\nonumber 
\end{eqnarray}
It easy to see that this process converges in $L^1$ to $Y_s$, uniformly in time, as $\var\rightarrow 0$. The $\F_{s(1-\var)}$-conditional density of the random variable $G^{\var}_s$ is normal with mean
\[E_G:=G_t+\int_t^{s(1-\var)}\mu_G(Y_{\l},\l)d\l\]
and covariance matrix
\[T_G(s,t):=(s-t)\s_G\s^{\top}_G>0.\]
Consider $\delta>0$, then there exists a countable collection of sets $Q_j\subset R^N\times[t^+,T]$ which, without loss of generality, can be supposed to be closed $N+1$-cubes, such that the projection of $Q_j$ over the coordinates $x_{m+1},\ldots,x_N$ lies in a bounded set of $A_2$ (this is not a restriction bacause we could always take the intersection of $\t{B}$ with an increasing sequence of sets whose projection over these coordinates is bounded), the intersection of a couple of cubes is a set of zero Lebesgue measure and
\[\t{B}\cap \left(\R^N\times[t^+,T]\right)\subset \cup_j Q_j\]
and
\[\sum_j {\cal L}^{N+1}(Q_j)<\delta.\]
We proceed at first with an estimate of the expected value (\ref{HAUS3.7}) where $\t{B}$ is substituted by $Q_j$ and $Y_s$ by $Y_s^{\var}$. By the distributional property of the approximating process we have
\begin{eqnarray}\label{HAUS4.7}
\lefteqn{\int_{t^+}^T\E_t\left[1_{Q_j}(Y^{\var}_s,s)\right]ds=}\nonumber\\
&=&\int_{t^+}^T\E_t\left[\E_{{s(1-\var)}}\left[1_{Q_j}(Y^{\var}_s,s)\right]\right]ds\nonumber\\
&=&\int_{t^+}^T\E_t\left[P\left((Y^{\var}_s,s)\in Q_j| \F_{s(1-\var)}\right)\right]ds.
\end{eqnarray}  
We shall use a diagonalization procedure to separate the variables in the following integration. Hence consider the orthogonal matrix $H$, such that $H^{\top}T_G H=Diag(\l_1,\ldots,\l_{N-m})$, where $\l_1,\ldots,\l_{N-m}$ are the eigenvalues of $T_G$.\\
If $Q_l(a)$, $B_{l}(a)$ respectively denote the  cube of centre $a$ and side of lenght $l>0$ and the open ball of radius $l$ and centre $a$ in $\R^{N-m}$, then again by the properties of $G_s^{\var}$, we have 
\begin{eqnarray}\label{HAUS5}
\lefteqn{P\left(G^{\var}_s\in Q_l(a)|\F_{s(1-\var)}\right)\leq}\nonumber\\
&\leq & \int_{B_{l\sqrt{N-m}}(a)}\frac{1}{\sqrt{2\pi}^{N-m} \sqrt{det(T_G(s,t))}}\exp\left(-\frac{1}{2}\la T^{-1}_G(s,t)(x-E_G),x-E_G\ra\right)dx\nonumber\\
&=&\int_{B_{l\sqrt{N-m}}(\o{a})}\frac{1}{\sqrt{2\pi(s-t)}^{N-m} \sqrt{\l_1\cdots\l_{N-m}}}\exp\left(-\frac{1}{2(s-t)}\sum_i \l_i^{-1}y_i^2\right)dy\nonumber\\
&\leq& \int_{Q_{l\sqrt{N-m}}(\o{a})}\frac{1}{\sqrt{2\pi(s-t)}^{N-m} \sqrt{\l_1\cdots\l_{N-m}}}\exp\left(-\frac{1}{2(s-t)}\sum_i \l_i^{-1}y_i^2\right)dy\nonumber\\
&=& \Pi_{i=1}^{N-m}\left[\Phi\left(\frac{l\sqrt{N-m}+\o{a}_i}{\sqrt{\l_i(s-t)}}\right)- \Phi\left(\frac{-l\sqrt{N-m}+\o{a}_i}{\sqrt{\l_i(s-t)}}\right)\right]\nonumber\\
&\leq& \frac{(2 l\sqrt{N-m})^{N-m}}{\sqrt{2\pi( t^+ -t)}\sqrt{det(\s_G\s_G^{\top})}}=C(N,m,t,t^+,\s_G){\cal L}^{N-m}(Q_l(a)),
\end{eqnarray}
where $\o{a}=H^{\top}(a-E_G)$, and $\Phi$ denotes the standard Gaussian distribution.\\
Writing  $Q_j=Q_j^1\times Q_j^2\times I_j$, for some cubes $Q_j^1\subset \R^m$, $Q_j^2\subset A_2$, and $I_j\subset[t^+,T]$, (\ref{HAUS4.7}) can be continued using the estimate (\ref{HAUS5}) in the following way
\begin{eqnarray}\label{HAUS6}
\lefteqn{=\int_{I_j}\E_t\left[ 1_{Q_j^1}(\pi^{M,\var}_s)P\left(G^{\var}_s\in Q_j^2|\F_{s(1-\var)}\right)\right]ds\leq}\nonumber\\
&\leq &C(N,m,t,t^+,\s_G) {\cal L}^{N-m}(Q_j^2)\int_{I_j}\E_t\left[ 1_{Q_j^1}(\pi^{M,\var}_s)\right]ds. 
\end{eqnarray}
By the condition (\ref{N>d}), there exists the density $d_{\pi}(x;\tau,\omega,t)\geq 0$ $P-a.s.$ of the measure $P(\pi_{\tau}^{-1}(\cdot)|X_t)$ with respect to ${\cal L}^m$ which, by the regularity of $X_t$, is a continuous function of $\tau\in [t^+,T]$.
By this property, (\ref{HAUS2}) and denoting $\widehat{M}$ the restriction of $M^{\top}$ over the kernel of $\s^{\top}$, we obtain
\begin{eqnarray}\label{HAUS7}
\lefteqn{\E_t\left[ 1_{Q_j^1}(\pi^{M,\var}_s)\right]=}\nonumber\\
&=& P\left( \pi^M_{s(1-\var)}\in Q_j^1|\F_{t}\right)\nonumber\\
&=& P\left(\pi_{s(1-\var)}\in \widehat{M}Q_j^1|\F_{t}\right)\nonumber\\
&=& \int_{\widehat{M}Q_j^1} d_{\pi}(x;s(1-\var),\omega,t)dx.
\end{eqnarray}
Introducing (\ref{HAUS7}) in (\ref{HAUS6}), we have
\begin{eqnarray}\label{HAUS8}
\lefteqn{\int_{t^+}^T\E_t\left[1_{Q_j}(Y^{\var}_s,s)\right]ds\leq}\nonumber\\ 
&= &  C(N,m,t,t^+,\s_G) {\cal L}^{N-m}(Q^2_j)\int_{I_j}\int_{\widehat{M}Q_j^1} d_{\pi}(x;s(1-\var),X_t,t)dx ds\nonumber\\
&= & C(N,m,t,t^+,\s_G)\int_{\widehat{M}Q_j^1\times Q_j^2\times I_j}1_{A_2}(y)d_{\pi}(x;s(1-\var),\omega,t)dx dy ds,\qquad\mbox{a.s.}
\end{eqnarray}
By the properties of the collection $\{Q_j\}_j$, we can write
\begin{eqnarray}\label{HAUS8.1}
\lefteqn{{\cal L}^{N+1}\left(\cup_j \widehat{M}Q_j^1\times Q_j^2\times I_j\right)=}\nonumber\\
&&|det\widehat{M}|{\cal L}^{N+1}(\cup_j Q_j)< |det\widehat{M}|\delta. 
\end{eqnarray}
For $P-a.s.$, the measure $\nu(\cdot)=\nu(\cdot;\omega,t)$, defined through the density function $(x,y,s)\mapsto 1_{A_2}(y)d_{\pi}(x;s(1-\var),\omega)$ over $\R^{N}\times[t^+,T]$, is an absolutely continuous measure with respect to $N+1$-dimensional Lebesgue measure. Therefore set $\beta>0$ and choose $\delta$ so that for every set $A$ with ${\cal L}^{N+1}(A)<\delta |det \widehat{M}|$, it holds $\nu(A)\leq\beta$. Hence using (\ref{HAUS8.1}), we deduce
\begin{eqnarray}\label{HAUS9}
\lefteqn{\E_t\left[\int_{t^+}^T 1_{\t{B}}(Y^{\var}_s,s)ds\right]\leq}\nonumber\\ &\leq& \nu\left(\cup_j \widehat{M}Q_j^1\times Q_j^2\times I_j\right)\leq \beta .
\end{eqnarray}
Then the assertion in the case $N>d$ follows by the arbitrary choice of $\beta$ and letting $\var\rightarrow 0$ and then $t^+\rightarrow t$ in (\ref{HAUS9}).\\
{\bf Case $N=d$.} Introducing the approximation
\begin{eqnarray}\label{HAUS10}
X^{\var}_s=X_t+\int_t^{s(1-\var)} \mu(X_{\l},T-\l)d\l+\s^{\var}(W_s-W_t),
\end{eqnarray}
where $\s^{\var}\rightarrow \s$, as $\var\rightarrow 0$ and $det(\s^{\var})\neq 0$, then we can repeat the same arguments used for the previous statement
substituting $X_s^{\var}$ at the process $Y_s^{\var}$ of the previous case.
Therefore we obtain
\begin{eqnarray}\label{HAUS11}
\E_t\left[\int_{t^+}^T 1_{B}(X^{\var}_s,s)ds\right]\leq C(\s^{\var})\delta,
\end{eqnarray}
for arbitrary $\delta>0$. Hence the limits with respect to $\var$ and $t^+$ is taken.
\cvd
\end{Dimo}

\begin{oss}\rm\label{osservazione}
The condition (\ref{N>d}) in some sense in optimally. Actually in the case $N>d$, the strong degeneration of the quadratic form $\s\s^{\top}$ could make false the property (\ref{HAUSSMANN0.1}). For istance if $N=2$ and $d=1$, choosing
\[\s^{\top}=(0,1),\;\;\;\mu=0,\;\;\; t=0,\;\;\;B=\{0\}\times (-\infty,0)\times [0,T],\]
then with $X_0=0$,
\[X_t=(0,W_t)\]
and 
\[\E\left[\int_0^T 1_B(X_s,s))ds\right]=\frac{T}{2}>0\]
which contradicts the assertion of Lemma \ref{HAUSSMANN0}.
\end{oss}

\begin{Dimo} {\bf Theorem \ref{HAUSSMANN}.}
Let $\Sigma(t)=\s(T-t)$, $\mu^{\circ}(x,t)=\mu(x,T-t)$, $h^{\circ}=h(x,T-t)$ and $\xi^{\circ}(t)=\xi(T-t)$, $r^{\circ}(t)=r(t)$ for every $(x,t)\in\R^N\times(0,T]$.
Then define the function $U^{\circ}(x,t)=U(x,T-t)$. $U^{\circ}$ is a viscosity solution of the equation
\begin{eqnarray}\label{DM2}
\lefteqn{\!\partial_t U^{\circ}\!+\!\frac{1}{2}tr(\Sigma\Sigma^{\top}\nabla^2 U^{\circ})\!+\!\la \mu^{\circ},\nabla U^{\circ} \ra\!-\!\rho\frac{|\Sigma^{\top}\nabla U^{\circ}|^2}{U^{\circ}+h^{\circ}+\xi^{\circ}(t)}=}\nonumber\\
&&=r^{\circ}(U^{\circ}\!+\!h^{\circ})\!-\!\tau h^{\circ}\qquad\mbox{and $U^{\circ}(x,T)\equiv0$.}
\end{eqnarray}
Moreover the equality (\ref{DM2}) holds for a.e. $(x,t)\in\R^N\times (0,T)$ with the weakly derivatives of $U^{\circ}$.
Applying Theorem 3.1 pg. 733, in \cite{ITO}, to the function $(x,s)\mapsto e^{-\int_t^s r^{\circ}(\kappa)d\kappa}U^{\circ}(x,s)$ relatively to the time interval $[t,T)$ there exist processes $\beta_0$, $\beta$, $\alpha$, such that
\begin{eqnarray}\label{ITO2}
0 &=& e^{-\int_t^T r^{\circ}(\kappa)d\kappa} U^{\circ}(X_T,T)=U^{\circ}(X_t,t)+\int_t^T\big[\beta_0(\omega,s)+\la\mu^{\circ}(X_s,s),\beta(\omega,s)\ra\nonumber\\
&&+\frac{1}{2}tr\big(\Sigma\Sigma^{\top}(s)\alpha(\omega,s)\big)\big]ds
+\int_t^T\beta^{\top}(\omega,s)\Sigma(s)\cdot dW_s,\qquad\mbox{$a.s.$}
\end{eqnarray}
Moreover there exists a set $A\subset\R^N\times(t,T]$ of full Lebesgue measure such that the usual derivatives of $U^{\circ}$
exist in $A$ and
\begin{eqnarray}\label{ITO3}
\beta_0(\omega,s)&=& e^{-\int_t^s r^{\circ}(\kappa)d\kappa}\left(\partial_s U^{\circ}(X_s(\omega),s)-r^{\circ}(s)U^{\circ}(X_s(\omega),s)\right),\nonumber\\
\beta(\omega,s)&=& e^{-\int_t^s r^{\circ}(\kappa)d\kappa}\nabla U^{\circ}(X_s(\omega),s),\\
\alpha(\omega,s)&=& e^{-\int_t^s r^{\circ}(\kappa)d\kappa}\nabla^2 U^{\circ}(X_s(\omega),s),\nonumber
\end{eqnarray}
whenever $(X_s(\omega),s)\in A$. Without loss of generality we can also assume that the equation (\ref{DM2}) holds over $A$.
The application of (\ref{DM2}), (\ref{ITO2}) and (\ref{ITO3}) yields 
\begin{eqnarray}\label{ITO4}
0 &=& U^{\circ}(X_t,t)+\int_t^T\big[\beta_0(s)+\la\mu^{\circ}(X_s,s),\beta(s)\ra+\frac{1}{2}tr\big(\Sigma\Sigma^{\top}(s)\alpha(s)\nonumber\\
&&-\beta^{\top}(s)\Sigma(s)\gamma_s\big)\big]1_{A}(X_s,s)ds+\int_t^T1_A(X_s,s)\beta^{\top}(\omega,s)\Sigma(s)\cdot d\widehat{W}_s+Z_{t,T}\nonumber\\
&=& U^{\circ}(X_t,t)+\int_t^T e^{-\int_t^s r^{\circ}(\kappa)d\kappa}\big[r^{\circ}(s)\big(U^{\circ}(X_s,s)+h^{\circ}(X_s,s)\big)-\tau h^{\circ}(X_s,s)\nonumber\\
&&-r^{\circ}(s)U^{\circ}(X_s,s)\big]ds+\int_t^T\beta^{\top}(\omega,s)\Sigma(s)1_{A}(X_s,s)\cdot d\widehat{W}_s\nonumber\\
&&+Z_{t,T},\qquad\qquad\qquad\mbox{a.s.}
\end{eqnarray}
By the boundness of the processes $\beta_0$, $\beta$, $\alpha$, we recognize that the remaining term $Z_{t,T}$ has the same structure of (\ref{HAUSSMANN0.1}) in Lemma \ref{HAUSSMANN0}, therefore, 
by the definition of $Q$, (\ref{ITO1}), as a measure equivalent to $P$ and Lemma \ref{HAUSSMANN0}, $Z_{t,T}$ has a null contidional $Q$-expected value.\\
Using this remark and taking the conditional expected value in both the left and right side of the equation (\ref{ITO4}), we obtain
\begin{eqnarray}\label{ITO5}
U^{\circ}(X_t,t)=\E_t^{Q}\left[\int_t^T e^{-\int_t^s r^{\circ}(\kappa)d\kappa}(\tau-r^{\circ}(s)) h^{\circ}(X_s,s)1_A(X_s,s)ds\right],\;\mbox{a.s.}
\end{eqnarray}
By the full measure of $A$ and noting that the functions $h$ and $r$ are defined everywhere and are in particular continuous, the assertion is proved.
\cvd
\end{Dimo}

Finally we obeserve that the same conclusions can be obtained for a process $X_t$ such that satisfies
\begin{eqnarray*}
dX_t=\mu(X_t,T-t) X_t dt+\s(T-t) X_t dW_t,\;\;X_i(0)>0,\;\i=1,\ldots,N.
\end{eqnarray*}
In this case can is possible to treat the problem using the same previous arguments simply making the change $y=log(x)$ in the equation (\ref{DM1}).\\ 
Also if the drift coefficient depends in a linear way by $x$, and the diffusion coefficient is constant with respect to $x$, then
we could repeat the same argument used before to state the regularity of $u$. Actually, in that case, we have a linear increasing with respect to the variable $x$, only in the first order linear term of the equation.
 

\end{document}